\documentclass[bj,12pt,preprint]{imsart}

\RequirePackage{amsthm,amsmath, enumerate}
\RequirePackage[colorlinks,citecolor=blue,urlcolor=blue]{hyperref}
\usepackage{mathrsfs}
\usepackage[margin=1in]{geometry}
\usepackage[english]{babel}
\usepackage{multicol}
\usepackage{pdfpages}
\usepackage{wrapfig}
\usepackage{tabularx}
\usepackage{caption}
\usepackage{bm}
\usepackage{mathrsfs}
\usepackage{amsmath}
\usepackage{enumitem}
\usepackage{pdfpages}
\usepackage{color}
\usepackage{float}
\usepackage{graphicx}
\usepackage{csquotes}
\usepackage{hyperref}
\usepackage{amsmath}

\usepackage{bbm}
\startlocaldefs
\numberwithin{equation}{section}
\theoremstyle{plain}
\newtheorem{thm}{Theorem}[section]
\theoremstyle{remark}
\newtheorem{rem}{Remark}
\theoremstyle{plain}

\theoremstyle{plain}
\newtheorem{lemma}{Lemma}
\theoremstyle{plain}
\newtheorem{proposition}{Proposition}
\theoremstyle{plain}

\theoremstyle{plain}
\newtheorem{example}{Example}
\theoremstyle{plain}

\endlocaldefs

\numberwithin{equation}{section}
\numberwithin{thm}{section}
\numberwithin{cor}{section}
\numberwithin{lemma}{section}
\numberwithin{proposition}{section}
\numberwithin{rem}{section}
\numberwithin{example}{section}
\numberwithin{res}{section}

\newcommand{\bl}{\textcolor{blue}}

\begin{document}

\begin{frontmatter}

\linespread{1}
\title{Second Order Correctness of Perturbation Bootstrap M-Estimator of 
Multiple Linear Regression Parameter\thanksref{T1}}
\runtitle{S.O.C. of Perturbation Bootstrap}
\thankstext{T1}{Research partially supported by NSF grants no. DMS 1310068, DMS 1613192}

\begin{aug}
\author{\fnms{Debraj} \snm{Das}$^{\text{a}}$\ead[label=c]{ddas25@wisc.edu}}
\and
\author{\fnms{S. N.} \snm{Lahiri}$^{\text{b}}$\ead[label=d]{snlahiri@ncsu.edu}}

\address{$^{\text{a}}$Department of Statistics, University of Wisconsin-Madison, 1300 University Avenue, Madison, WI 53706,
USA.
\printead{c}}
\address{$^{\text{b}}$Department of Statistics, North Carolina State University, 2311 Stinson Dr, Raleigh, NC 27695-8203,
USA.
\printead{d}}

\runauthor{Das D. and Lahiri S. N.}

\end{aug}

\linespread{1.25}

\begin{abstract}
. Consider the multiple linear regression model 
$y_{i} = \mathbf{x'_{i}} \bm{\beta} + \epsilon_{i}$, 
where $\epsilon_i$'s are independent and identically distributed random variables, $\mathbf{x_i}$'s are known design vectors and $\bm{\beta}$ is the $p \times 1$ vector of parameters. An effective way of approximating the distribution of the M-estimator $\bm{\bar{\beta}_n}$, after proper centering and scaling, is the Perturbation Bootstrap Method. In this current work, second order results of this non-naive bootstrap method have been investigated. Second order correctness  is important for reducing the approximation error uniformly to $o(n^{-1/2})$ to get better inferences. We show that the classical studentized version of the bootstrapped estimator fails to be second order correct. We introduce an innovative modification in the studentized version of the bootstrapped statistic and show that the modified bootstrapped pivot is second order correct (S.O.C.) for approximating the distribution of the studentized M-estimator. Additionally, we show that the Perturbation Bootstrap continues to be S.O.C. when the errors $\epsilon_i$'s are independent, but may not be identically distributed. These findings establish perturbation Bootstrap approximation as a significant improvement over asymptotic normality in the regression M-estimation.
\end{abstract}

\begin{keyword}
\kwd{M-Estimation, S.O.C., Perturbation Bootstrap, Edgeworth Expansion, Studentization, Residual Bootstrap, Generalized Bootstrap, Wild Bootstrap}
\end{keyword}

\end{frontmatter}

\section{ Introduction}
%


Consider the multiple linear regression model :

\begin{equation}
y_{i} = \mathbf{x'_{i}}\bm{\beta} + \epsilon_{i}, \; \;\;\;\;     i = 1,2,\ldots,n
\end{equation}
where $y_1,\ldots,y_n$ are responses, $\epsilon_1,\ldots,\epsilon_n$ are independent and identically distributed (IID) random variables with common distribution $F$ (say), $\mathbf{x_1},\ldots,\mathbf{x_n}$ are known non random design vectors and $\bm{\beta}$ is the $p$-dimensional vector of parameters.

\par
Suppose $\bm{\bar{\beta}}_n$ is the M-estimator of $\bm{\beta}$ corresponding to the objective function $\Lambda(\cdot)$ i.e. $\bm{\bar{\beta}_n} = \operatorname*{arg\,min}_{\mathbf{t}} \sum_{i=1}^{n}$ $\Lambda(y_i - \mathbf{x'_i t})$. Now if $\psi(\cdot)$ is the derivative of $\Lambda(\cdot)$, then $\bm{\bar{\beta}_n}$ is the M-estimator corresponding to the score function $\psi(\cdot)$ and is defined as the solution of the vector equation

\begin{equation*}
\sum_{i = 1}^{n}\mathbf{x_i}\psi(y_i - \mathbf{x'_i}\bm{\beta}) = \mathbf{0}.
\end{equation*}
 It is known [cf. Huber(1981)] that under some conditions on the objective function, design vectors and error distribution $F$; $(\bm{\bar{\beta}_n} - \bm{\beta})$ with proper scaling has an asymptotically normal distribution with mean $\mathbf{0}$ and dispersion matrix $\sigma^2 \mathbf{I_p}$ where $\sigma^2 = \mathbf{E}\psi^2(\epsilon_1)/{\mathbf{E}^2\psi'(\epsilon_1)}$. 

\par
After introduction of bootstrap by Efron in 1979 as a resampling technique, it has been widely used as a distributional approximation method. Resampling from the naive empirical distribution of the centered residuals in a regression setup, called residual bootstrap, was introduced by Freedman (1981).  Freedman (1981) and Bickel and Freedman (1981b) had shown that given data, the conditional distribution of $\sqrt{n}(\bm{\beta^*_n}-\bm{\bar{\beta}_n})$ converges to the same normal distribution as the distribution of $\sqrt{n}(\bm{\bar{\beta}_n}-\bm{\beta})$ when $\bm{\bar{\beta}_n}$ is the usual least square estimator of $\bm{\beta}$, that is, when $\Lambda(x)=x^2$. It implies that the residual bootstrap approximation to the exact distribution of the least square estimator is first order correct as in the case of normal approximation. The advantage of the residual bootstrap approximation over normal approximation for the distribution of linear contrasts of least square estimator for general $p$ was first shown by Navidi (1989) by investigating the underlying Edgeworth Expansion (EE); although heuristics behind the same was given by Liu (1988) in restricted case $p=1$. Consequently, EE for the general M-estimator of $\bm{\beta}$ was obtained by Lahiri (1989b) when $p=1$; whereas the same for the multivariate least square estimator was found by Qumsiyeh (1990a). EE of standardized and studentized versions of the general M-estimator in multiple linear regression setup was first obtained by Lahiri (1992). Lahiri (1992) also established the second order results for residual bootstrap in regression M-estimation.

A natural generalization of sampling from the naive empirical distribution is to sample from a weighted empirical distribution to obtain the bootstrap sample residuals. Broadly, the resulting bootstrap procedure is called the weighted or generalized bootstrap. It was introduced by Mason and Newton (1992) for bootstrapping mean of a collection of IID random variables. Mason and Newton (1992) considered exchangeable weights and established its consistency. Lahiri (1992) established second order correctness of generalized bootstrap in approximating the distribution of the M-estimator for the model (1.1) when the weights are chosen in a particular fashion depending on the design vectors. Wellner and Zhan (1996) proved the consistency of infinite dimensional generalized bootstrapped M-estimators. Consequently, Chatterjee and Bose (2005) established distributional consistency of generalized bootstrap in estimating equations and showed that generalized bootstrap can be used in order to estimate the asymptotic variance of the original estimator. Chatterjee and Bose (2005) also mentioned the bias correction essential for achieving second order correctness. An important special case of generalized bootstrap is the bayesian bootstrap of Rubin (1981). Rao and Zhao (1992) showed that the distribution function of M-estimator for the model (1.1) can be approximated consistently by bayesian bootstrap. See the monograph of Barbe and Bertail (2012) for an extensive study of generalized bootstrap.

A close relative to the generalized bootstrap procedure is the wild bootstrap. It was introduced by Wu (1986) in multiple linear regression model (1.1) with errors $\epsilon_i$'s being heteroscedastic. Beran (1986) justified wild bootstrap method by pointing out that the distribution of the least square estimator can be approximated consistently by the wild bootstrap approximation. Second order results of wild bootstrap in heteroscedastic regression model was first established by Liu (1988) when $p=1$. Liu (1988) also showed that usual residual bootstrap is not capable of approximating the distribution of the least square estimator upto second order in heteroscedastic setup and described a modification in resampling procedure which can establish second order correctness. For general $p$, the heuristics behind achieving second order correctness by wild bootstrap in homoscedastic least square regression were discussed in Mammen (1993). Recently, Kline and Santos (2011) developed a score based bootstrap method depending on wild bootstrap in M-estimation for the homoscedastic model (1.1) and established consistency of the procedure for Wald and Lagrange Multiplier type tests for a class of M-estimators under misspecification and clustering of data.

A novel bootstrap technique, called the perturbation bootstrap was introduced by Jin, Ying, and Wei (2001) as a resampling procedure where the objective function having a U-process structure was perturbed by non-negative random quantities. Jin, Ying, and Wei (2001) showed that in standardized setup, the conditional distribution of the perturbation resampling estimator given the data and the distribution of the original estimator have the same limiting distribution which means this resampling method is first order correct without studentization. In a recent work, Minnier, Tian, and Cai (2011) also applied this perturbation resampling method in penalized regression setup such as Adaptive Lasso, SCAD, $l_q$ penalty and showed that the standardized perturbed penalized estimator is first order correct. But, second order properties of this new bootstrap method have remained largely unexplored in the context of multiple linear regression. In this current work, the perturbation bootstrap approximation is shown to be S.O.C. for the distribution of studentized M-estimator for the regression model (1.1). An extension to the case of independent and non-IID errors is also established, showing the robustness of perturbation bootstrap towards the presence of heteroscedasticity. Therefore, besides the existing bootstrap methods, the perturbation bootstrap method can also be used in regression M-estimation for making inferences regarding the regression parameters and higher order accuracy can be achieved than the normal approximation. 

A classical way of studentization in bootstrap setup, in case of regression M-estimator and for IID errors, is to consider the studentization factor to be $\sigma_n^* = s_n^*\tau_n^{*-1}$, $\tau_n^*=n^{-1}\sum_{i=1}^{n}\psi'(\epsilon_i^*)$, $s^{*2}_n = n^{-1}\sum_{i=1}^{n}\psi^2(\epsilon_i^*)$ where $\epsilon_i^*=y_i-\mathbf{x_i'}\bm{\beta_n^*}$, $i\in\{1,\dots,n\}$, with $\bm{\beta_n^*}$ being the perturbation bootstrapped estimator of $\bm{\beta}$, defined in Section 2. Although the residual bootstrapped estimator is S.O.C. after straight-forward studentization, the same pivot fails to be S.O.C. in the case of perturbation bootstrap. Two important special cases are considered as examples in this respect. The reason behind this failure is that although the bootstrap residuals are sufficient in capturing the variability of the bootstrapped estimator in residual bootstrap, it is not enough in the case of perturbation resampling. Modifications have been proposed as remedies and are shown to be S.O.C. The modifications are based on the novel idea that the variability of the random perturbing quantities $G_i^*$ ($1\leq i\leq n$) along with the bootstrap residuals are required to capture the variability of the perturbation bootstrapped estimator; whereas individually they are not sufficient. For technical details, see Section 4.2 and Section 5.1.

\par
With a view to establish second order correctness, we start with the standardized setup and then proceed to studentization. First, we find a two-term EE of the conditional density of a suitable stochastic approximation of the concerned bootstrapped pivot and then we show that it is the required two-term EE corresponding to the bootstrapped pivot. The result then follows by comparing the EE of the bootstrapped pivot with that of underlying original pivot. The techniques that are to be used in finding EE have been demonstrated and discussed in Bhattacharya and Ghosh  (1978), Bhattacharya and Rao (1986), Navidi (1989) and Lahiri (1992).

A significant volume of work is available in bootstrapping M-estimators. We will conclude this section by briefly reviewing the literature. Bootstrapping M-estimators in linear model has been studied by Navidi(1989), Lahiri (1992, 1996), Rao and Zhao (1992), Qumsiyeh (1994), Karabulut and Lahiri (1997), Jin, Ying and Wei (2001), Hu (2001), El Bantli (2004) among others. And in the applications other than linear model, bootstrapping in M-estimation and its subclasses has been investigated by Arcones and Gin\'e (1992), Lahiri (1994), Wellner and Zhan (1996), Allen and Datta (1999), Hu and Kalbfleisch (2000), Hlavka (2003), Wang and Zhou (2004), Chatterjee and Bose (2005), Ma and Kosorok (2005), Lahiri and Zhu (2006), Cheng and Huang (2010), Feng et. al. (2011), Lee (2012), Cheng (2015), among others.

The rest of the paper is organized as follows. Perturbation bootstrap is described briefly in Section 2. Section 3 states the assumptions and motivations behind considering those assumptions. Main results for IID case, along with the modification in bootstrap studentization, are stated in Section 4. An extension to the case of independent and non-IID errors is proposed in Section 5. An outline of the proofs are given in Section 6. Section 7 states concluding remarks. The details of the proofs are available in a supplementary material Das and Lahiri (2017).

\section{Description of Perturbation Bootstrap}

In the perturbation bootstrap, the objective function $\Lambda(\cdot)$ has been perturbed several times by a non-negative random quantity to get a bootstrapped estimate of $\bm{\beta}$. It has nothing to do with residuals in resampling stage, unlike the residual and weighted bootstrap. More precisely, the perturbation bootstrap estimator $\bm{\beta_n^*}$ is defined as 
\begin{equation*}
\bm{\beta^*_n} = \operatorname*{arg\,min}_{\mathbf{t}} \sum_{i=1}^{n}\Lambda(y_i - \mathbf{x'_i t})G^*_i
\end{equation*}
or in terms of the score function $\psi(\cdot)$, as the solution of the vector equation
\begin{equation}
\sum_{i = 1}^{n}\mathbf{x_i}\psi(y_i - \mathbf{x'_i}\bm{\beta})G^*_i = \mathbf{0}
\end{equation}
where $G^*_i, i \in \{1, \ldots, n\}$ are non-negative and non-degenerate completely known random variables, considered as perturbation quantities. Note that, if $\mu_{G^*}$ is the mean of $G_1^*$, then $\bm{\bar{\beta}_n}$ is the solution of $\mathbf{E}\Big(\sum_{i = 1}^{n}\mathbf{x_i}\psi(\bar{\epsilon}_i)$ $G^*_i|\epsilon_1,\dots, \epsilon_n\Big)=\sum_{i = 1}^{n}\mathbf{x_i}\psi(\bar{\epsilon}_i)\mu_{G^*}=0$ where $\bar{\epsilon}_i=y_i-\mathbf{x_i'}\bar{\bm{\beta}}_n$, $i\in\{1,\dots,n\}$, are the residuals corresponding to the M-estimator $\bar{\bm{\beta}}_n$. This observation will be helpful in finding a suitable stochastic approximation in bootstrap regime. For details, see Section 6.

The central idea of the perturbation bootstrap is to draw a relatively large collection of IID random samples $\{(G^{*b}_{1},\ldots,G^{*b}_{n}): b = 1,\ldots,B\}$ from the distribution of $G^*_1$ and then to find the conditional empirical distribution of $\sqrt{n}(\bm{\beta^*_n}-\bm{\bar{\beta}_n)}$ given data ${y_i: i = 1,\ldots,n}$, by solving 
\begin{equation*}
\sum_{i = 1}^{n}\mathbf{x_i}\psi(y_i - \mathbf{x'_i}\bm{\beta})G^{*b}_i = \mathbf{0}
\end{equation*}
for each $b\in \{1,\ldots,B\}$; to approximate the distribution of $\sqrt{n}(\bm{\bar{\beta}_n}-\bm{\beta})$ asymptotically. As a result the bootstrapped distribution may be used as an approximation to the original distribution, just like the normal approximation, in constructing confidence intervals and testing of hypotheses regarding $\bm{\beta}$.

Now, in the perturbation bootstrap M-estimation, $G_i^*$'s can be thought of as weight corresponding to the $i$th data point $(\mathbf{x_i},y_i)$. To make it easier to understand, consider the least square setup i.e. $\Lambda(x) = x^2$. In this case $\bm{\beta^*_n}$ takes the form
\begin{equation}
\bm{\beta^*_n} = \Big(\sum_{i=1}^{n}\mathbf{x_ix'_i}G^*_i\Big)^{-1}\Big(\sum_{i=1}^{n}\mathbf{x_i}y_iG^*_i\Big)
\end{equation}
indicating that the perturbing quantities $G^*_i$'s can be thought of as weights.

\begin{rem}
Consider the least square estimator $\bm{\hat{\beta}_n}$. Then keeping the asymptotic properties fixed, the perturbation bootstrap version $\bm{\hat{\beta}_{1n}}^*$ of $\bm{\hat{\beta}_n}$ can be defined alternatively as the solution of 
\begin{equation*}
\sum_{i = 1}^{n}\mathbf{x_i}(y_i - \mathbf{x'_i}\bm{\beta})\big(G^*_i-\mu_{G^*}\big) + \sum_{i = 1}^{n}\mathbf{x_i}\mathbf{x'_i}(\bm{\hat{\beta}_n} - \bm{\beta})\big(2\mu_{G^*}-G_i^*\big) = \mathbf{0}
\end{equation*}
which in turn implies that $\bm{\hat{\beta}_{1n}}^*$ is the solution of
\begin{equation}
\sum_{i = 1}^{n}\mathbf{x_i}(z_i^* - \mathbf{x'_i}\bm{\beta}) = \mathbf{0}
\end{equation}
where $z_i^*=\mathbf{x'_i}\bm{\hat{\beta}_n}+\hat{\epsilon}_i[\mu_{G^*}^{-1}(G_i^*-\mu_{G^*})]$, $\hat{\epsilon}_i=y_i-\mathbf{x'_i}\bm{\hat{\beta}_n}$, $i\in \{1,\dots,n\}$. On the other hand, the simple wild bootstrap version $\bm{\hat{\beta}_{2n}}^*$ of $\bm{\hat{\beta}_n}$ is defined as the solution of 
\begin{equation}
\sum_{i = 1}^{n}\mathbf{x_i}(y_i^* - \mathbf{x'_i}\bm{\beta}) = \mathbf{0}
\end{equation}
where $y_i^*=\mathbf{x'_i}\bm{\hat{\beta}_n}+\hat{\epsilon}_it_i$, $i\in\{1,\dots,n\}$ and $\{t_1,\dots,t_n\}$ is a set of IID random variables independent of $\{\epsilon_1,\dots,\epsilon_n\}$ with $\mathbf{E}t_1=0$, $\mathbf{Var}(t_1)=1$. Additionally, one needs $\mathbf{E}(t_1^3)=1$ for establishing second order correctness of wild bootstrap approximation [cf. Liu (1988), Mammen (1993)]. Now  Looking at (2.3) and (2.4) and in view of assumption (A.5)(ii), it can be said that the perturbation bootstrap coincides with the wild bootstrap in least square setup. Therefore
one can view perturbation bootstrap as a generalization of the wild bootstrap in regression M-estimation.
\end{rem}

\begin{rem}
There is a basic difference between perturbation bootstrap and weighted bootstrap with respect to the construction of the bootstrapped estimator. Whereas in the perturbation bootstrap, the bootstrapped estimator is defined through the non-negative and non-degenerate random perturbations of the objective function; in weighted bootstrap, the bootstrapped estimator is defined through bootstrap samples drawn from a weighted empirical distribution. See for example the construction of the weighted bootstrapped estimator corresponding to Theorem 2.3 of Lahiri (1992) and compare it with our construction as stated in Section 2. However, as pointed out by a referee, one can think of the perturbation bootstrap, defined in Section 2, as the weighted bootstrap version of some statistical functional if the design vectors are random. Suppose, $\{(x_1,y_1)\dots,(x_n,y_n)\}$ are IID with underlying probability measure $\mathbf{Q}$. Then one can write
\begin{align*}
\bm{\beta}=T(\mathbf{Q}) = \operatorname*{arg\,min}_{\bm{t}} \mathbf{E}_{\mathbf{Q}}\Big{[}\Lambda(y_i - \bm{x'_i t})\Big{]}
\end{align*}
for some statistical functional $T(\cdot)$. Define empirical measures  $\mathbf{Q}_{n}=n^{-1}\sum_{i=1}^{n}\mathbbm{1}{(x_i,y_i)}$ and $\mathbf{Q}_{n,W}$ $=n^{-1}\sum_{i=1}^{n}\mathbbm{1}{(x_i,y_i)}W_i$ where $\mathbbm{1}(\cdot)$ is the indicator function and $\{W_1,\dots,W_n\}$ are weights. Then we have $\bar{\bm{\beta}}_n=T(\mathbf{Q}_n)$ and $\bm{\beta}_n^*=T(\mathbf{Q}_{n,W})$ when $W_i=G_i^*$, $i\in \{1,\dots,n\}$. The weighted bootstrap of general statistical functionals of only the IID random variables is considered in the monograph of Barbe and Bertail (2012). Second order correctness of weighted bootstrap of  standardized mean of IID random variables was established by Haeusler et. al. (1992) under two choices of weights. One choice is the non-negative IID weights and the other one is the self-normalized sum of non-negative IID random variables. Their results were extended by Barbe and Bertail (2012) for general statistical functionals in IID case when the weights are self-normalized sum of non-negative IID random variables [cf. Corollary 4.1 of Barbe and Bertail (2012)]. For general M-estimation, Chatterjee (1999) showed that weighted bootstrap estimator is generally biased and established its second order correctness after properly correcting for the bias.
 To the best of our knowledge, there is no second order result available in the literature under studentized setup for general statistical functional.
 In this article, we have assumed the design vectors to be non-random, implying that our setup fits neither in the general statistical functional setup of Barbe and Bertail (2012) nor in the general M-estimation setup of Chatterjee (1999); although Theorem 5.1 continue to hold when the design is random. Throughout the article we consider weights to be non-negative IID. Our main motivation is to explore second order results in studentized setup which, unlike the standardized (i.e., the known variance) case,  is applicable in practice. Further, we prove our results in  the situation when errors are heteroscedastic. We establish all our second order correctness results without requiring any bias correction.

\end{rem}

\section{Assumptions}
Suppose, $\mathbf{x_i} = (x_{i1},x_{i2},\ldots,x_{ip})'$. Define, $\mathbf{D_n} \equiv \mathbf{D} = (\sum_{i = 1}^{n}\mathbf{x_i}\mathbf{x'_i})^{1/2}$, $\mathbf{A_n} = n^{-1}\mathbf{D^{2}}$, $\mathbf{d_i} = \mathbf{D^{-1}x_i}$, $1\leq i \leq n$ and $q = \dfrac{p(p+1)}{2}$. Also define, $q\times 1$ vector $\mathbf{z_i} = (x_{i1}^2,x_{i1}x_{i2},\ldots,x_{i1}x_{ip}$ $,x_{i2}^2,x_{i2}x_{i3},\ldots$ $, x_{i2}x_{ip},\ldots,x_{ip}^2)'$ .
Note that for any constants $a_i,\ldots,a_n \in \mathscr{R}$, $\sum_{i = 1}^{n}a_i\mathbf{z_i} = \mathbf{0}$ which implies and is implied by $\sum_{i=1}^{n}a_i\mathbf{x_ix'_i} = \mathbf{0}$. Hence, $\{\mathbf{z_1},\ldots,\mathbf{z_n}\}$ are linearly independent if and only if $\{\mathbf{x_ix'_i} : 1\leq i\leq n\}$ are linearly independent. Therefore, $r_n =$ the rank of $\sum_{i = 1}^{n}\mathbf{z_iz'_i}$ is nondecreasing in n. So, if $r = \max\{r_n :n\geq 1\}$ then without loss of generality (w.l.g.), we can assume that $r_n = r$ for all $n\geq q$. Consider canonical decomposition of $\sum_{i = 1}^{n}\mathbf{z_iz'_i}$ as

\begin{equation*}
\mathbf{L}(\sum_{i = 1}^{n}\mathbf{z_i z'_i})\mathbf{L'} = \begin{bmatrix}
\mathbf{I_r} \;\;\;\mathbf{0}\\
\mathbf{0}\;\;\;\; \mathbf{0}
\end{bmatrix}
\end{equation*}
where $\mathbf{L}$ is a $q\times q$ non-singular matrix. Partition $\mathbf{L}$ as $\mathbf{L'} = [\mathbf{L'_1}\;\; \mathbf{L'_2}]$, where $\mathbf{L_1}$ is of order $r \times q$. Define $r \times 1$ vector $\mathbf{\tilde{z}_i}$ by 
\begin{equation*}
\mathbf{\tilde{z}_i} = \mathbf{L_1z_i},\;\;\;\;\; 1\leq i\leq n
\end{equation*}
Note that $\sum_{i=1}^{n}\mathbf{\tilde{z}_i\tilde{z}'_i} = \mathbf{L_1}(\sum_{i=1}^{n}\mathbf{z_iz'_i})\mathbf{L'_1} = \mathbf{I_r}$. Suppose, $\mathbf{v_i}=\mathbf{(x'_i\psi(\epsilon_1),z'_i\psi'(\epsilon_1))'}$. $\mathbf{\breve{z}_i} =(\mathbf{z'_i},n^{-1})'$.

\par
Let, $\mathbf{\Phi_V}$ denotes the normal distribution with mean $\mathbf{0}$ and dispersion matrix $\mathbf{V}$ and $\mathbf{\phi_V}$ is the density of $\mathbf{\Phi_V}$. Write $\mathbf{\Phi_V }= \mathbf{\Phi}$ and $\mathbf{\phi_V} = \mathbf{\phi}$ when $\mathbf{V}$ is the identity matrix. $h',h''$ denote respectively first and second derivatives of real valued function $h$ that is twice differentiable. Also $||.||$ denotes euclidean norm.For any set $B\in \mathscr{R}^p$ and $\epsilon>0$, $\delta B$ denotes the boundary of $B$, $|B|$ denotes the cardinality of $B$ and $B^{\epsilon}=\{\mathbf{x}: \mathbf{x}\in \mathscr{R}^p \;\text{and}\; d(\mathbf{x},B)<\epsilon\}$ where $d(\mathbf{x},B)=\inf\{||\mathbf{x}-\mathbf{y}||:\mathbf{y}\in B\}$. For a function $f :\ \mathscr{R}^l\ \rightarrow \ \mathscr{R}$ and a non-negative integral vector $\bm{\alpha} = (\alpha_1, \alpha_2,\ldots,\alpha_l)'$, $D^{\bm{\alpha}}f = D_1^{\alpha_1}\ldots D_l^{\alpha_l}f$, where $D_j^{\alpha_j}f$ denotes $\alpha_j$ times partial derivative of $f$ with respect to the $j$th component of its argument, $1\leq j \leq l$. Also assume that $(\mathbf{e_1},\ldots, \mathbf{e_p})'$ is the standard basis of $\mathscr{R}^p$. Let, $\mathbf{P_*}$ and $\mathbf{E_*}$ respectively denote conditional bootstrap probability and conditional expectation of $G_1^{*}$ given data. The class of sets $\mathscr {B}$ denotes the collection of borel subsets of $\mathscr{R}^p$ satisfying 
\begin{equation}
\sup\limits_{B \in  \mathscr{B}}\; \mathbf{\Phi}((\delta B)^{\epsilon}) = O(\epsilon) \;\;\; as \;\epsilon \downarrow 0
\end{equation}  Next we state the assumptions:

\vspace{2mm}
\begin{enumerate}[label=(A.\arabic*)]
\item $\psi(\cdot)$ is twice differentiable and $\psi^{\prime\prime}(\cdot)$ satisfies a Lipschitz condition of order $\alpha$ for some $0< 2\alpha\leq 1$.
\item \begin{enumerate}[label=(\roman*)] \item $\mathbf{A_n} \rightarrow \mathbf{A_1}$ as $ n\rightarrow \infty$ for some positive definite matrix $\mathbf{A_1}$.
\item  $\mathbf{E}(n^{-1}\sum_{i=1}^{n}\mathbf{v_iv'_i}) \rightarrow \mathbf{A_2}$  as $ n\rightarrow \infty$
for some non-singular 
matrix 
$\mathbf{A_2}$, where expectation is with respect to  $F$.\\
\hspace*{-8mm}(ii)$'$ $\mathbf{E}(n^{-1}\sum_{i=1}^{n}\mathbf{\tilde{v}_i\tilde{v}_i'}) \rightarrow \mathbf{A_3}$  as $ n\rightarrow \infty$
for some non-singular matrix 
$\mathbf{A_3}$ where $\mathbf{\tilde{v}_i}$ is defined as same way as $\mathbf{v_i}$ with $\mathbf{z_i}$ being replaced by $\mathbf{\breve{z}_i}$.
\item $n^{\alpha/2}(\sum_{i = 1}^{n}||\mathbf{d_i}||^{6 +2\alpha})^{1/2} + \sum_{i=1}^{n}||\mathbf{\tilde{z}_i}||^4= O(n^{-1})$
\end{enumerate}
\item  \begin{enumerate}[label=(\roman*)]
\item $\mathbf{E}\psi(\epsilon_1) = 0$ and $\sigma^2 = \mathbf{E}\psi^2(\epsilon_1)/\mathbf{E}(\psi'(\epsilon_1)) \in (0,\infty)$.
\item $\mathbf{E}|\psi(\epsilon_1)|^4 + \mathbf{E}|\psi'(\epsilon_1)|^4 + \mathbf{E}|\psi''(\epsilon_1)|^2 < \infty$.
\end{enumerate}
\item $G_{i}^{*}$ and $\epsilon_i$ are independent for all $1\leq i\leq n$.
\item \begin{enumerate}[label=(\roman*)] \item $\mathbf{E}G_{1}^{*3} < \infty$ 
\item $\mathbf{Var}(G_1^*) = \mu_{G^*}^2$, $\mathbf{E}(G_1^* - \mu_{G^*})^3 = \mu_{G^*}^3$.
\item $\big(G^*_1-\mu_{G^*}\big)$ satisfies Cramer's condition:\\ 
\hspace*{10mm}$\limsup_{|t|\rightarrow \infty}\big|\mathbf{E}\big(exp\big(it\big(G^*_1-\mu_{G^*}\big)\big)\big)\big|<1$.\\ 
\hspace*{-9mm}(iii)$'$ $\big(\big(G^*_1-\mu_{G^*}\big), \big(G_1^* - \mu_{G^*}\big)^{2}\big)$ satisfies Cramer's condition:\\ 
\hspace*{5mm}$\limsup_{||(t_1,t_2)||\rightarrow \infty}\Big|\mathbf{E}\big(exp\big(it_1\big(G^*_1-\mu_{G^*}\big)+it_2\big(G_1^* - \mu_{G^*}\big)^{2}\big)\big)\Big|<1$ 
\end{enumerate}
\item 
\begin{enumerate}[label=(\roman*)]
\item $\big(\psi(\epsilon_1), \psi'(\epsilon_1)\big)$ satisfies Cramer's condition:\\ 
\hspace*{5mm}$\limsup_{||(t_1,t_2)||\rightarrow \infty}\Big|\mathbf{E}\big(exp\big(it_1\psi(\epsilon_1)+it_2\psi'(\epsilon_1)\big)\big)\Big|<1$ \\ 
\hspace*{-7mm}(i)$'$ $\big(\psi(\epsilon_1), \psi'(\epsilon_1), \psi^2(\epsilon_1)\big)$ satisfies Cramer's condition:\\ 
\hspace*{1mm}$\limsup_{||(t_1,t_2,t_3)||\rightarrow \infty}\Big|\mathbf{E}\big(exp\big(it_1\psi(\epsilon_1)+it_2\psi'(\epsilon_1)+it_3\psi^2(\epsilon_1)\big)\big)\Big|<1$ 
\end{enumerate}
\end{enumerate}
\hspace{3mm}

Define $\mathbf{\bar{v}_i}=(\mathbf{\bar{x}'}_i,\mathbf{\bar{z}'_i})'$ where $\mathbf{\bar{x}_i}=\mathbf{x_i}\psi(\bar{\epsilon}_i)$, $\mathbf{\bar{z}_i}=\mathbf{z_i}\psi'(\bar{\epsilon}_i)$; $\{\bar{\epsilon}_1,\ldots,\bar{\epsilon}_n\}$ being the set of residuals. Also, define $\mathbf{\bar{A}_{2n}}=n^{-1}\sum_{i=1}^{n}\mathbf{\bar{x}_i\bar{x}'_i}$ and $\mathbf{\bar{A}_{1n}}=n^{-1}\sum_{i=1}^{n}\mathbf{x_ix'_i}\psi'(\bar{\epsilon}_i)$. Note that   $n^{-1}\sum_{i=1}^{n}\mathbf{\bar{v}_i\bar{v}'_i}$ is an estimate of the matrix $\mathbf{E}(n^{-1}\sum_{i=1}^{n}\mathbf{v_iv'_i})$ and due to assumption (A.2)(ii), $\sum_{i=1}^{n}\mathbf{\bar{v}_i\bar{v}'_i}$ is non-singular for sufficiently large $n$. Hence, without loss of generality the canonical decomposition of $\sum_{i = 1}^{n}\mathbf{\bar{v}_i\bar{v}'_i}$ can be assumed as
\begin{equation*}
\mathbf{B}\Big(\sum_{i = 1}^{n}\mathbf{\bar{v}_i\bar{v}'_i}\Big)\mathbf{B'} = \mathbf{I_k}
\end{equation*}
where $k=p+q$ and $\mathbf{B}$ is a $k\times k$ non-singular matrix. Define $k \times 1$ vector $\mathbf{\breve{v}_i}$ by 
\begin{equation*}
\mathbf{\breve{v}_i} = \mathbf{B\bar{v}_i},\;\;\;\;\; 1\leq i\leq n
\end{equation*}

To find valid EE in the perturbation bootstrap regime, the following condition [cf. Navidi (1989)] is also required:
\vspace*{2mm}

\hspace*{-4mm}(A.7) There exists a $\delta > 0$ such that $-K_n(\delta)/log\gamma_n \rightarrow \infty$
where $\mathbf{B_n(\delta)} = \{1\leq i \leq n : \hspace*{9mm}  (\mathbf{\breve{v}'_it})^2 > \delta\gamma_n^2$ for all $\mathbf{t}\in \mathscr{R}^{k}$ with $||\mathbf{t}||^2 = 1\}$, $K_n(\delta) = |\mathbf{B_n(\delta)}|$, the cardinality of \hspace*{9mm}  the set $\mathbf{B_n(\delta)}$, and
$\gamma_n = (\sum_{i=1}^{n}||\mathbf{\breve{v}_i}||^4)^{1/2}$.\\

But note that the condition (A.7) has already been satisfied in our set up due to Lemma 6.2 and the proposition in Lahiri (1992).\\

Now we briefly explain the assumptions. Assumption (A.1) is smoothness condition on the score function $\psi(\cdot)$. This condition is essential for obtaining a Taylor's expansion of $\psi(\cdot)$ around regression errors. Assumption (A.2) presents the regularity conditions on the design vectors necessary to find EE. For the validity of asymptotic normality of the regression M-estimator, only (A.2)(i) is enough [cf. Huber (1981)]; whereas additional condition (A.2)(ii) is required for the validity of the EE. (A.2)(iii) states atmost how fast the $L^2$ norm of the design vectors can increase to get a valid EE. This condition is somewhat stronger than the condition (C.6) assumed in Lahiri (1992); although there was a reduction in accuracy of bootstrap approximation due to this relaxation. This type of conditions are quite common in the literature of edgeworth expansions in regression setup; see for example Navidi (1989), Qumsiyeh (1990a). We now state an example where assumption (A.2) (iii) is fulfilled. 
\begin{example}
Suppose, $\{\mathbf{X^{(1)}},\dots, \mathbf{X^{(p)}}\}$ is a set of independent random vectors where $\mathbf{X}^{(j)}=(X_{1j},\dots,$ $X_{nj})^\prime$ is a vector of $n$ IID copies of the non-degenerate random variable $X_{1j}$, $j \in \{1,\dots,p\}$. Define, $p\times p$ matrix $\mathbf{M}=((m_{jk}))_{j,k=1,\dots,p}$ where $m_{jk}=\mathbf{E}(X_{1j}^2X_{1k}^2)$ and $n\times p$ matrix $\mathbf{X}=\big{(}\mathbf{X^{(1)}},\dots, \mathbf{X^{(p)}}\big{)}$. Assume, $\mathbf{E}(X_{1j})=\mathbf{E}(X_{1j}^3)=0$ and $\mathbf{E}|X_{1j}|^{\bl{8}}<\infty$ for all $j\in \{1,\dots,p\}$ and $\det(M)\neq 0$. Then for the design matrix $\mathbf{X}$, assumption \emph{(A.2) (iii)} holds with probability 1 (w.p. 1).
\end{example}

\subsubsection*{ proof :} 
For the design matrix $\mathbf{X}$, $\mathbf{x_i} = (X_{i1},X_{i2},\ldots,X_{ip})'$ and $\mathbf{z_i} = (X_{i1}^2,X_{i1}X_{i2},\ldots,X_{i1}X_{ip},X_{i2}^2$ $,X_{i2}X_{i3},$ $\ldots,$ $ X_{i2}X_{ip},\ldots,X_{ip}^2)'$ for $i\in \{1,\dots,n\}$.

First note that if all the entries of $\mathbf{X}$ are IID then the condition $\det(M)\neq 0$ is redundant. By Kolmogorov strong law of large numbers, $\mathbf{A_n}=n^{-1}\mathbf{D^2}$ $\rightarrow diag\big{(}$ $\mathbf{E}(X_{11}^2),\dots,$ $\mathbf{E}(X_{1p}^2)\big{)}$ and $n^{-1}\sum_{i=1}^{n}||\mathbf{x_{i}}||^{6+2\alpha}\rightarrow \mathbf{E}||\mathbf{x_{1}}||^{6+2\alpha}$ both $w.p. 1$ and hence
\begin{align}
n^{\alpha/2}\Big(\sum_{i = 1}^{n}||\mathbf{d_i}||^{6 +2\alpha}\Big)^{1/2}& \leq n^{\alpha/2}||\mathbf{D^{-1}}||^{3+\alpha}\Big{(}\sum_{i=1}^{n}||\mathbf{x_{i}}||^{6+2\alpha}\Big{)}^{1/2} \nonumber \\
& =O(n^{-1})\;\;\; w.p. 1
\end{align}

Again, since $\mathbf{M}$ is a non-singular matrix, $n^{-1}\sum_{i=1}^{n}\mathbf{z}_i \mathbf{z}_i^\prime\rightarrow \mathbf{N}\; w.p. 1,$ for some positive definite matrix $\mathbf{N}$. This implies that $||\mathbf{L}||=O(n^{-1/2})\; w.p. 1$ and hence
\begin{align}
\sum_{i=1}^{n}||\mathbf{\tilde{z}_i}||^4& \leq ||\mathbf{L}||^4 \sum_{i=1}^{n}||\mathbf{z_i}||^4 \nonumber \\
& =O(n^{-1})\;\;\; w.p. 1
\end{align}
Therefore, our claim follows from (3.2) and (3.3).\\

Assumption (A.3) is the moment condition on the error variables through the score function $\psi(\cdot)$. (A.3)(i) is generally assumed to establish asymptotic normality. Assumption (A.4) is inherent in the present setup, since $G_i^*$'s are introduced by us to define the bootstrapped estimator whereas $\epsilon_i$'s are already present in the process of data generation. The conditions present in Assumption (A.5) are moment and smoothness conditions on the perturbing quantities $G_i^*$'s, required for the valid two term EE in bootstrap setup.  The Cramer's condition is very common in the literature of edgeworth expansions. Cramer's condition is satisfied when the distribution of $(G_1^*-\mu_{G^*})$ or $((G_1^*-\mu_{G}^*), (G_1^*-EG_1^*)^2)$ has a non-degenerate component which is absolutely continuous with respect to Lebesgue measure [cf. Hall (1992)]. An immediate choice of the distribution of $G_1^*$ is $Beta$($\gamma$, $\delta$) where $3\gamma=\delta=3/2$. Also one can investigate $Generalized \;Beta$ family of distributions for more choices of the distribution of $G_1^*$. Assumption (A.6) is the Cramer's condition on the errors. Although this assumption is not needed for obtaining EE of the bootstrapped estimators, it is needed for obtaining EE for the original M-estimator. 

Note that the condition (A.7) is somewhat abstract. Hence as pointed out by a referee, some clarification would be helpful. To this end, it is worth mentioning that to find formal EE for the standardized bootstrapped pivot (see section 4.1), the most difficult step is to show
\begin{align}
\max_{|\bm{\alpha}|\leq p+q+4}{\int_{C_1\leq \gamma_n||\mathbf{t}||\leq C_2}|D^{\bm{\alpha}}\mathbf{E_*}e^{i\mathbf{t}^\prime\mathbf{T}_n^*}|d\bm{t}}=o_p\big(n^{-1/2}\big)
\end{align}
where $C_1, C_2$ are non-negative constants and $\mathbf{T}_n^*=\sum_{i=1}^{n}\big{(}\mathbf{\breve{X}_i^*}-\mathbf{E_*}(\mathbf{\breve{X}_i^*})\big{)}$, with $\mathbf{\breve{X}_i^*}=\mathbf{\breve{v}_i}(G_i-\mu_{G^*})\mathbf{1}\big{(}||\mathbf{\breve{v}_i}(G_i-\mu_{G^*})||\leq 1\big{)}$. Now it is easy to see that for any $|\bm{\alpha}|\leq p+q+4$,  $|D^{\bm{\alpha}}\mathbf{E_*}e^{i\mathbf{t}^\prime\mathbf{T}_n^*}|$ is bounded above by a sum of $n^{|\alpha|}$-terms, each of which is bounded above by
\begin{align*}
C(\alpha) \cdot \max\{\mathbf{E_*}||\mathbf{\breve{X}_i^*}-\mathbf{E_*}(\mathbf{\breve{X}_i^*})||^{|\bm{\alpha}|}: i \in \mathbf{I_n^*}\} \cdot \prod_{i \in \mathbf{I^{*c}_n}}|\mathbf{E_*}e^{i\mathbf{t}^\prime\mathbf{\breve{X}_i^*}}|
\end{align*}
where $\mathbf{I_n^*}\subset \{1,\dots, n\}$ is of size $|\bm{\alpha}|$ and $\mathbf{I^{*c}_n}=\{1,\dots,n\}\backslash \mathbf{I_n^*}$ and $C(\bm{\alpha})$ is a constant which depends only on $\bm{\alpha}$.

Now note that for all $i \in \{1,\dots,n\}$,
\begin{align*}
&\mathbf{E_*}||\mathbf{\breve{X}_i^*}-\mathbf{E_*}(\mathbf{\breve{X}_i^*})||^{|\bm{\alpha}|}\leq 2^{|\bm{\alpha}|}\\
\text {and}\;\;\;\;\;\; & |\mathbf{E_*}e^{i\mathbf{t}^\prime\mathbf{\breve{X}_i^*}}|\leq |\mathbf{E_*}e^{i\mathbf{t}^\prime\mathbf{\breve{v}_i}(G_i-\mu_{G^*})}| + 2\mathbf{P_*}\big{(}||\mathbf{\breve{v}_i}(G_i-\mu_{G^*})||> 1\big{)}
\end{align*}
Hence, in view of Cramer's condition (A.5) (iii) and Lemma 6.2, if there exists a sequence of sets $\{\mathbf{J_n}\}_{n\geq 1}$ such that $\mathbf{J_n}\subset \{1,\dots,n\}$ and for all $i\in \mathbf{J_n}$, $\gamma_n^{-1}|\mathbf{t}^\prime\mathbf{\breve{v}_i}|> \xi$ for some $\xi>0$, then for some $0<\theta<1$ we have
\begin{align}
&\sup\Big{\{}\prod_{i \in \mathbf{I^{*c}_n}}|\mathbf{E_*}e^{i\mathbf{t}^\prime\mathbf{\breve{X}_i^*}}|:C_1\leq \gamma_n||\mathbf{t}||\leq C_2\Big{\}} \nonumber\\
&\leq \sup\Big{\{}\prod_{i \in \mathbf{I^{*c}_n}\cap \mathbf{J_n}}|\mathbf{E_*}e^{i\gamma_n^{-1}\mathbf{t}^\prime\mathbf{\breve{X}_i^*}}|:C_1\leq ||\mathbf{t}||\leq C_2\Big{\}} \nonumber\\
& \leq \theta^{|\mathbf{I^{*c}_n}\cap \mathbf{J_n}|}
\end{align}
Again $|\mathbf{I^{*c}_n}\cap \mathbf{J_n}|\geq |\mathbf{J_n}|-|\bm{\alpha}|$ and $\gamma_n\geq kn^{-1}$. Therefore, to achieve (3.4), it is enough to have
\begin{align*}
n^{2(p+q)+4}\cdot \theta^{|\mathbf{J_n}|-(p+q+4)} = o(n^{-1/2})
\end{align*}
Hence due to Lemma 6.2, it is enough to have $|\mathbf{J_n}|\geq a_n - C\cdot \log \gamma_n$ for some positive constant $C$ and a sequence of constants $\{a_n\}$ increasing to $\infty$. This observation together with (3.5) justifies condition (A.7).\\    

We will denote the assumptions (A.1)-(A.5) by (A.1)$'$-(A.5)$'$ when (A.2) and (A.5) are respectively defined with (ii)$'$ and (iii)$'$ instead of (ii) and (iii).

\section{Main Results}

\subsection{Rate of Perturbation Bootstrap Approximation}
Here we will state the approximation results both in standardized and studentized setup. It is well known that $\sqrt{n}\bm{\bar{\beta}_n}$ has asymptotic variance $\sigma^2\mathbf{A_n^{-1}}$. So, the standardized version of the M-estimator $\bm{\bar{\beta}_n}$ is defined as $\mathbf{F_n} = \sqrt{n}\sigma^{-1}\mathbf{A_n^{1/2}}(\bm{\bar{\beta}_n} - \bm{\beta})$. Now to define the standardized version of the corresponding bootstrapped statistic $\bm{\beta_n^*}$, we need its conditional asymptotic variance, given the data. Using Taylor's expansion, it is quite easy to get the conditional asymptotic variance of $\sqrt{n}\bm{\beta_n^*}$ as $\mathbf{\bar{A}_{1n}^{-1}\bar{A}_{2n}\bar{A}_{1n}^{-1}}$. Note that inverse of the matrices $\mathbf{\bar{A}_{1n}^{-1}}$ and $\mathbf{\bar{A}_{2n}^{-1}}$ are well defined for sufficiently large sample size n due to the assumption (A.2)(i) and (A.3)(ii). Hence, the standardized bootstrapped M-estimator $\mathbf{F_n^*}$ can be defined as
\begin{equation*}
\mathbf{F_n^*}=\sqrt{n}\mathbf{\bar{\Sigma}_n^{-1/2}}(\bm{\beta^*_n}-\bm{\bar{\beta}_n})
\end{equation*}
where $\mathbf{\bar{\Sigma}_n^{-1/2}}=\mathbf{\bar{A}_{2n}^{-1/2}\bar{A}_{1n}}$, $\mathbf{\bar{A}_{2n}^{1/2}}$ being defined in terms of the spectral decomposition of $\mathbf{\bar{A}_{2n}}$; although it can be defined in many different ways [cf. Lahiri (1994)]. Under some regularity conditions, both the distribution of  $\mathbf{F_n}$ and the conditional distribution of $\mathbf{F_n^*}$ can be shown to be approximated asymptotically by a Normal distribution with mean $\mathbf{0}$ and variance $\mathbf{I_p}$. Hence, it is straightforward that perturbation bootstrap approximation to the distribution of the M-estimator is first order correct. The second order result in standardized case is formally stated in Theorem 4.1.

\begin{proposition}
Suppose, the assumptions \emph{(A.1)-(A4), (A.5)(i)} hold. Then there exist constant $C_1>0$ and a sequence of Borel sets $\mathbf{Q_{1n}}\subseteq \mathscr{R}^n$, such that $\mathbf{P}((\epsilon_1,\ldots,\epsilon_n)\in \mathbf{Q_{1n}})\rightarrow 1$ as $n\rightarrow\infty$, and given $(\epsilon_1,\ldots,\epsilon_n)\in \mathbf{Q_{1n}}$, $n\geq C_1$ such that there exists a sequence of statistics $\{\bm{\beta_n^*}\}_{n\geq 1}$ such that
\begin{equation*}
\mathbf{P_*}\big(\bm{\beta_n^*} \; solves  \;(2.1)\; and \;||\bm{\beta_n^*} - \bm{\bar{\beta}_n}||\leq C_1.n^{-1/2}.(log n)^{1/2}\big) \geq 1 - \delta_nn^{-1/2}
\end{equation*} 
where $\delta_n \equiv \delta_n(\epsilon_1,\ldots,\epsilon_n)$ tends to 0.
\end{proposition}

\begin{thm}
Let $\{\bm{\beta_n^*}\}_{n\geq 1}$ be a sequence of statistics satisfying Proposition $\text{4.1}$ depending on $(\epsilon_1,\ldots,\epsilon_n)$. Assume, the assumptions \emph{(A.1)-(A.5)} hold. 
\begin{itemize}
\item[\emph{(a)}] Then there exist constant $C_2>0$ and a sequence of Borel sets $\mathbf{Q_{2n}}\subseteq \mathscr{R}^n$ and polynomial $a^*_n(\cdot,\psi,G^{*})$ depending on first three moments of $G_1^{*}$ and on $\psi(\cdot)$, $\psi^\prime(\cdot)$ \& $\psi^{\prime\prime}(\cdot)$ through the residuals $\{\bar{\epsilon}_1,\dots,\bar{\epsilon}_n\}$ such that given $(\epsilon_1,\ldots,\epsilon_n)\in \mathbf{Q_{2n}}$, with $\mathbf{P}((\epsilon_1,\ldots
,\epsilon_n)\in \mathbf{Q_{2n}})\rightarrow 1$, we have for $n\geq C_2$, 
\begin{equation*}
\sup\limits_{B \in  \mathscr{B}} |\mathbf{P_*}(\mathbf{F_n^*}\in B) - \int_B\xi^*_n(\mathbf{x})d\mathbf{x}| \leq  \delta_n n^{-1/2}
\end{equation*}
where $\xi^*_n(\mathbf{x}) = (1 + n^{-1/2}a^*_n(\mathbf{x},\psi,G^{*}))\phi(\mathbf{x})$ and $\delta_n \equiv \delta_n(\epsilon_1,\ldots,\epsilon_n)$ tends to 0.\\

\item[\emph{(b)}] Suppose in addition assumption \emph{(A.6)(i)} holds. Then we have,
\begin{equation*}
\sup\limits_{B \in  \mathscr{B}} \big{|}\mathbf{P_*}(\mathbf{F_n^*}\in B) - \mathbf{P}(\mathbf{F_n}\in B)\big{|} = o_p( n^{-1/2})
\end{equation*}
\end{itemize}
\end{thm}

Now, the quantity $\sigma^2$ is mostly unavailable in practical circumstances. Hence, the non-pivotal quantity like $\mathbf{F_n}$ is very rare in use in providing valid inferences. It is more reasonable to explore the asymptotic properties of a pivotal quantity, like the studentized version of the M-estimator $\bm{\bar{\beta}_n}$. Depending on the observed residuals $\bar{\epsilon}_i=y_i - \mathbf{x_i}'\bm{\bar{\beta}_n},\;i\in\{1,\ldots,n\}$, the natural way to define an estimator of $\sigma^2$ is $\hat{\sigma}_n^{2}$ where  $\hat{\sigma}_n= s_n\tau_n^{-1}$,  $\tau_n=n^{-1}\sum_{i=1}^{n}\psi'(\bar{\epsilon}_i)$ and $s^{2}_n = n^{-1}\sum_{i=1}^{n}\psi^2(\bar{\epsilon}_i)$. Hence, the studentized M-estimator in regression setup may be defined as $\mathbf{H_n}=\sqrt{n}\hat{\sigma}_n^{-1}\mathbf{A_n^{1/2}}(\bm{\bar{\beta}_n}-\bm{\beta})$. Define the studentized version of the corresponding bootstrapped estimator as
\begin{equation*}
\mathbf{H_n^*}=\sqrt{n}\sigma_n^{*-1}\hat{\sigma}_n\mathbf{\bar{\Sigma}_n^{-1/2}}(\bm{\beta_n^*} - \bm{\bar{\beta}_n})
\end{equation*}
where $\sigma_n^* = s_n^*\tau_n^{*-1}$, $\tau_n^*=n^{-1}\sum_{i=1}^{n}\psi'(\epsilon_i^*)$, $s^{*2}_n = n^{-1}\sum_{i=1}^{n}\psi^2(\epsilon_i^*)$ and $\hat{\sigma}_n^2$ and $\mathbf{\bar{\Sigma}_n^{-1/2}}$ are as defined earlier.\\

\begin{thm}
Suppose, the assumptions \emph{(A.1)-(A.5)} hold. 
\begin{itemize}
\item[\emph{(a)}]Then there exist constant $C_3>0$ and a sequence of Borel sets $\mathbf{Q_{3n}}\subseteq \mathscr{R}^n$ and polynomial $\tilde{a}^*_n(\cdot,\psi,G^{*})$ depending on first three moments of $G_1^{*}$ and on $\psi(\cdot)$, $\psi^\prime(\cdot)$ \& $\psi^{\prime\prime}(\cdot)$ through the residuals $\{\bar{\epsilon}_1,\dots,\bar{\epsilon}_n\}$, such that given $(\epsilon_1,\ldots,\epsilon_n)\in \mathbf{Q_{3n}}$, with $\mathbf{P}((\epsilon_1,\ldots
,\epsilon_n)\in \mathbf{Q_{3n}})\rightarrow 1$, we have for $n\geq C_3$, 
\begin{equation*}
\sup\limits_{B \in  \mathscr{B}} |\mathbf{P_*}(\mathbf{H_n^*}\in B) - \int_B\tilde{\xi}^*_n(\mathbf{x})d\mathbf{x}| \leq \delta_n n^{-1/2}
\end{equation*}
where $\tilde{\xi}^*_n(\mathbf{x}) = (1 + n^{-1/2}\tilde{a}^*_n(\mathbf{x},\psi,G^{*}))\phi(\mathbf{x})$ and $\delta_n \equiv \delta_n(\epsilon_1,\ldots,\epsilon_n)$ tends to 0.\\

\hspace{-8mm}Suppose in addition assumption \emph{(A.6)(i)}$'$ holds. Then\\
\item[\emph{(b)}]
 for the collection of Borel sets $\mathscr {B}$ defined by \emph{(3.1)},
\begin{equation*}
\sup\limits_{B \in  \mathscr{B}} \big{|}\mathbf{P_*}(\mathbf{H_n^*}\in B) - \mathbf{P}(\mathbf{H_n}\in B)\big{|} = O_p( n^{-1/2})
\end{equation*}

\item[\emph{(c)}]
if $\;2E\psi^2(\epsilon_1)E\psi(\epsilon_1)\psi'(\epsilon_1)\neq E\psi'(\epsilon_1)E\psi^3(\epsilon_1)$, then there exists $\epsilon>0$ such that,
\begin{align*}
\mathbf{P}\Big(\liminf_{n\rightarrow\infty}\sqrt{n}\Big[\sup\limits_{B \in  \mathscr{B}} \big{|}\mathbf{P_*}(\mathbf{H_n^*}\in B) - \mathbf{P}(\mathbf{H_n}\in B)\big{|}\Big] > \epsilon\Big) =1
\end{align*}
\end{itemize}

\end{thm}

\begin{rem}
 Proposition 4.1 states that there exists a sequence of perturbation bootstrapped estimator $\bm{\beta^*_n}$ within a neighborhood of length $C.n^{-1/2}(logn)^{1/2}$ around the original M-estimator $\bm{\bar{\beta}_n}$ outside a set of bootstrap probability $o_p(n^{-1/2})$. This existence result is essential in finding valid EEs in bootstrap regime. This can be compared with Theorem 2.3 (a) of Lahiri (1992), where similar kind of result was shown in case of residual and generalized bootstrap.
\end{rem}

\begin{rem}
 Note that, where as the error term in approximating the distribution of M-estimator by perturbation bootstrap is of order $O_p(n^{-1/2})$ in the prevalent studentize setup, it reduces the order of the error of approximation to $o_p(n^{-1/2})$ in simple standardized setup. This means that the difference between coefficients corresponding to the term $n^{-1/2}$ in the EEs of original and bootstrapped estimator can be made arbitrarily small in standardized setup, but not in usual studentized setup.
\end{rem}

\begin{rem}
 To understand part (c) of Theorem 4.2, consider the usual least square estimator. In least square setup, the condition in the Theorem 4.2 (c) reduces to $E\epsilon^3 \neq 0$. This simply means that if the studentization in perturbation bootstrapped version is performed analogously as in case of original least square estimator, then the bootstrap distribution can not correct the original distribution upto second order. If this is investigated more deeply, then it can be observed that the usual studentized perturbation bootstrap approximation can not correct for the skewness of the error distribution $F$.
\end{rem}

\subsubsection{\bf{Examples}}
Theorem 4.2 concludes that the standard way of performing studentization of the bootstrapped estimator is first order correct. In order to show that the usual studentized setup is not second order correct, we consider following two important special cases with $\psi(x)=x$.  

\subsection*{{Example 4.1}}
Consider the observations $\{y_1,\ldots,y_n\}$ are coming from the distribution $F$ with a location shift $\mu$. This in terms of regression model becomes
\begin{equation*}
y_i = \mu + \epsilon_i
\end{equation*}
Hence, in this setup $p=1$, $\bm{\beta} = \mu$ and $x_i = 1$ for all $i\in\{1,\ldots,n\}$.\\ 

It can be shown that in this setup, $\tilde{\xi}_n(\cdot)$ and $\tilde{\xi}_n^*(\cdot)$, the EE of $\mathbf{H_n}$ and $\mathbf{H_n^*}$ respectively, turn out to be
\begin{equation*}
\tilde{\xi}_n(x) = \Big[1 - n^{-1/2}\big{\{}\tilde{b}_{11}\dfrac{d}{dx}+6^{-1}\tilde{b}_{31}\dfrac{d^3}{dx^3}\big{\}}\Big]\phi(x)
\end{equation*}
\begin{equation*}
\tilde{\xi}_n^*(x) = \Big[1 - n^{-1/2}\big{\{}\tilde{b}_{11}^*\dfrac{d}{dx}+6^{-1}\tilde{b}_{31}^*\dfrac{d^3}{dx^3}\big{\}}\Big]\phi(x)
\end{equation*}
where\\

$\tilde{b}_{11}=-2^{-1}\sigma^{-3}E\epsilon_1^3$, $\tilde{b}_{31}=-2\sigma^{-3}E\epsilon_1^3$

$\tilde{b}_{11}^*=-2\sigma_n^{-1}n^{-1}\sum_{i=1}^{n}\bar{\epsilon}_i$, $\tilde{b}_{31}^*=\sigma_n^{-3}-12\sigma_n^{-1}n^{-1}\sum_{i=1}^{n}\bar{\epsilon}_i$\\

It is clear that $\tilde{b}_{11}^*$ as well as $\tilde{b}_{31}^*$ are not converging respectively to $\tilde{b}_{11}$ and $\tilde{b}_{31}$ in probability and hence the perturbation bootstrap method is not second order correct in the above setup when the bootstrapped estimator is studentized in the usual manner.

\subsection*{{Example 4.2}}
Consider the simple linear regression model
\begin{align*}
y_i = \beta_0 + \beta_1x_i + \epsilon_i
\end{align*}
where $\beta_0$ and $\beta_1$ are parameters of interest and $\epsilon_i$'s are IID errors. This model, in terms of our multivariate linear regression structure, can be written as $y_i=\mathbf{\tilde{x}_i}'\bm{\beta}+\epsilon_i$ where $\bm{\beta}=(\beta_0,\beta_1)'$ and $\mathbf{\tilde{x}_i}=(1,x_i)'$. Hence, the  EEs of the original and bootstrapped estimators upto the order $o(n^{-1/2})$, after usual studentization, respectively becomes
\begin{equation*}
\tilde{\xi}_n(y_1,y_2) = \Bigg[1 - n^{-1/2}\Big{\{}\sum_{j=1}^{2}\tilde{b}_{11}^{*(j)}\dfrac{\partial}{\partial y_j}+\sum_{j=0}^{3}\dfrac{\tilde{b}_{31}^{(j,3-j)}}{j!(3-j)!}D^{(j,3-j)}\Big{\}}\Bigg]\phi(y_1,y_2)
\end{equation*}
\begin{equation*}
\tilde{\xi}_n^*(y_1,y_2) = \Bigg[1 - n^{-1/2}\Big{\{}\sum_{j=1}^{2}\tilde{b}_{11}^{*(j)}\dfrac{\partial}{\partial y_j}+\sum_{j=0}^{3}\dfrac{\tilde{b}_{31}^{*(j,3-j)}}{j!(3-j)!}D^{(j,3-j)}\Big{\}}\Bigg]\phi(y_1,y_2)
\end{equation*}
where\\

$\tilde{b}_{11}^{(j)}=-2^{-1}\Big[ n^{-1}\sum_{i=1}^{n}\mathbf{e'_jA_n^{-1/2}\tilde{x}_i}\Big]\gamma_1$\\

$\tilde{b}_{11}^{*(j)}=o_p(1)$\\

where $\big{(}\mathbf{e_1},\dots,\mathbf{e_p}\big{)}^\prime$ is the standard basis of $\mathscr{R}^p$, $j=1 \;\text{or}\; 2$, $\gamma_1$ is the coefficient of skewness of $\epsilon_1$,
$\mathbf{A_n} = n^{-1}\sum_{i=1}^{n}\mathbf{\tilde{x}_i\tilde{x}'_i}=\begin{bmatrix}
1 \;\;\;\;\;\bar{x}\\
\bar{x}\;\;\;\; \bar{x^2}
\end{bmatrix}$
where $\bar{x}=n^{-1}\sum_{i=1}^{n}x_i$ and $\bar{x^2}=n^{-1}\sum_{i=1}^{n}x_i^2$. $\mathbf{\bar{A}_{2n}}$ is as defined in general setup with $\mathbf{\tilde{x}_i}$ in place of $x_i$ for $i\in \{1,\dots,n\}$. The form of the coefficients $\tilde{b}_{31}^{(j_1,j_2)}$ and $\tilde{b}_{31}^{*(j_1,j_2)}$ are given in the supplementary material Das and Lahiri (2017) for all $(j_1,j_2)\in\{(a,b):a,b\in\{0,1,2,3\}\; \text{and}\; a+b=3\}$.

Note that, the coefficients $\tilde{b}_{11}^{(j)}$, $1\leq j \leq p$, all can not vanish together unless $\gamma_1=0$ and hence $\tilde{b}_{11}^{*(j)}$ can not converge to $\tilde{b}_{11}^{(j)}$ unless $\gamma_1=0$. Similarly, it can be shown that same condition is required to have the closeness of the coefficients $\tilde{b}_{31}^{(j,3-j)}$ and $\tilde{b}_{31}^{*(j,3-j)}$. Hence, the two EEs can not get closer unless $\gamma_1=0$, similar to the Example 4.1. This is exactly what is stated in the part (c) of Theorem 4.2 in most general form.

\subsection{ Modification to the bootstrapped pivot}
As it has been seen that $\mathbf{H_n^*}$, the usual studentized version of the perturbation bootstrapped estimator is not attending the desired optimal rate $o_p(n^{-1/2})$, so in the perspective of statistical inference, perturbation bootstrap is not advantageous over asymptotic normal approximation. For the sake of obtaining second order correctness, define the modified studentized $\bm{\beta^*_n}$ as
\begin{equation}
\mathbf{\tilde{H}_n^*} = \sqrt{n}(\tilde{\sigma}_n^*)^{-1}\hat{\sigma}_n\mathbf{\bar{\Sigma}_n^{-1/2}}(\bm{\beta^*_n}-\bm{\bar{\beta}_n})
\end{equation}
where\\ 

$\tilde{\sigma}_n^* = \tilde{s}_n^*\tilde{\tau}_n^{*-1}$, $\tilde{\tau}_n^*=n^{-1}\sum_{i=1}^{n}\psi'(\epsilon_i^*)G_i^*$, $\tilde{s}^{*2}_n = n^{-1}\sum_{i=1}^{n}\psi^2(\epsilon_i^*)(G_i^*-\mu_{G^*})^2$.\\

The bootstrapped statistic $\mathbf{\tilde{H}_n^*}$ can be seen to be achieving the optimal rate, namely $o_p(n^{-1/2})$, in approximating the original studentized M-estimator $\mathbf{H_n}$, which is formally stated in the following theorem:

\begin{thm}
Suppose, the assumptions \emph{(A.1)}$'$\emph{-(A.5)}$'$ hold. Also assume $EG_1^{*4}< \infty$.
 \begin{itemize}
 \item[\emph{(a)}] Then there exist constant $C_4>0$ and a sequence of Borel sets $\mathbf{Q_{4n}}\subseteq \mathscr{R}^n$ and polynomial $\bar{a}^*_n(\cdot,\psi,G^{*})$ depending on first three moments of $G_1^{*}$ and on $\psi(\cdot)$, $\psi^\prime(\cdot)$ \& $\psi^{\prime\prime}(\cdot)$ through the residuals $\{\bar{\epsilon}_1,\dots,\bar{\epsilon}_n\}$, such that given $(\epsilon_1,\ldots,\epsilon_n)\in \mathbf{Q_{4n}}$, with $\mathbf{P}((\epsilon_1,\ldots
,\epsilon_n)\in \mathbf{Q_{4n}})\rightarrow 1$, we have for $n\geq C_4$, 
\begin{equation*}
\sup\limits_{B \in \mathscr{B}} |\mathbf{P_*}(\mathbf{\tilde{H}_n^*}\in B) - \int_B\bar{\xi}^*_n(\mathbf{x})d\mathbf{x}| \leq \delta_n n^{-1/2}
\end{equation*}
where $\bar{\xi}^*_n(\mathbf{x}) = (1 + n^{-1/2}\bar{a}^*_n(\mathbf{x},\psi,G^{*}))\phi(\mathbf{x})$ and $\delta_n \equiv \delta_n(\epsilon_1,\ldots,\epsilon_n)$ tends to 0.\\

\item[\emph{(b)}]
Suppose, in addition \emph{(A.6)(i)}$'$ holds. Then, for the collection of Borel sets defined by \emph{(3.1)},
\begin{equation*}
\sup\limits_{B \in  \mathscr{B}} \big{|}\mathbf{P_*}(\mathbf{\tilde{H}_n^*}\in B) - \mathbf{P}(\mathbf{H_n}\in B)\big{|} = o_p(n^{-1/2})
\end{equation*}
\end{itemize}

\end{thm}

\begin{rem}
The modification that is needed to make the perturbation bootstrap method correct upto second order, suggests that besides incorporating the effect of bootstrap randomization through $\psi(\cdot)$ and $\psi'(\cdot)$ in the studentization factor of the bootstrap estimator, it is also essential to blend properly the effect of randomization that is coming directly from the perturbing quantities $G_i^*$s. 
\end{rem}

\begin{rem}
As pointed out by a referee, the usefulness of the above results depend critically on the rate of the probability $\mathbf{P}\big{(}(\epsilon_1,\dots,\epsilon_n \in \mathbf{Q_{in}})\big{)}$, $i=1,2,3,4$. Following the steps of the proofs, it can be shown that $\mathbf{P}\big{(}(\epsilon_1,\dots,\epsilon_n \in \mathbf{Q_{n}})\big{)}=1-O\big{(}n^{-1/2}(\log n)^{-2+\gamma_2}\big{)}$ where $\mathbf{Q_n}=\cap_{i=1}^{4}\mathbf{Q_{in}}$, for some $\gamma_2 \in (0,2)$, although the rate can be improved under moment condition stronger than (A.3) (ii). In general, if $\mathbf{E}|\psi(\epsilon_1)|^{2\gamma_3} + \mathbf{E}|\psi'(\epsilon_1)|^{2\gamma_3} + \mathbf{E}|\psi''(\epsilon_1)|^{\gamma_3} < \infty$ for some natural number $\gamma_3 \geq 2$, then analogously it can be shown that $\mathbf{P}\big{(}(\epsilon_1,\dots,\epsilon_n \in \mathbf{Q_{n}})\big{)}=1-O\big{(}n^{-(2\gamma_3-3)/2}(\log n)^{-\gamma_3+\gamma_2}\big{)}$ for some $\gamma_2 \in (0,\gamma_3)$. This implies that second order correctness of perturbation bootstrap can be established in almost sure sense under higher moment condition. 
\end{rem}

\begin{rem}
The condition (3.1) on the collection of Borel subsets $\mathscr{B}$ of $\mathscr{R}^p$, that is considered in the above theorems, is somewhat abstract. This condition is needed for achieving two goals. One is to obtain valid EE for the normalized part of the underlying pivot [cf. Corollary 20.4 of Bhattacharya and Rao (1986)] and the other one is to bound the remainder term with an order $o(n^{-1/2})$ with probability (or bootstrap probability) $1-o(n^{-1/2})$. These two together allow us to get EE for the underlying pivots. A natural choice for $\mathscr{B}$ is the collection of all Borel measurable convex subsets of $\mathscr{R}^p$.   
\end{rem}

\section{Extension to independent and non-identically distributed errors}
In this section, we will extend second order results of perturbation bootstrap to the model (1.1) with independent and non-identically distributed [hereafter referred to as non-IID] errors. Clearly the case of non-IID errors includes the situation when the regression errors are heteroscedastic. In many practical situations, the measurements obtained have different variability due to a number of reasons and hence it is crucial for an inference procedure to be robust towards the presence of heteroscedasticity. We will show that perturbation bootstrap can approximate the exact distribution of the regression M-estimator $\bm{\bar{\beta}_n}$ upto second order even when the errors are non-IID.

Before stating second order result in non-IID case, we describe briefly the literature available on bootstrap methods in heteroscedastic regression. Although there is huge literature available on bootstrap in homoscedastic regression, literature on bootstrap in heteroscedastic regression models is limited. Wu (1986) mentioned the limitation of residual bootstrap in heteroscedasticity and introduced wild bootstrap in least square regression. Beran (1986) gave justification behind consistency of wild bootstrap. Liu (1988) established second order correctness of wild bootstrap in heteroscedastic least square regression when dimension $p=1$. Liu (1988) proposed a modification of residual bootstrap in resampling stage and gave justification behind second order correctness. You and Chen (2006) proved consistency of wild bootstrap in approximating the distribution of least square estimator in semiparametric heteroscedastic regression model. Davidson and Flachaire (2008) and Davidson and Mackkinnon (2010) developed  wild bootstrap procedure for testing the coefficients in heteroscedastic linear regression. Arlot (2009) developed a resampling-based penalization procedure for model selection based on exchangeable weighted bootstrap.

\vspace*{1mm}

We state some additional assumptions needed to establish second order correctness. Define, $\mathbf{A_{1n}}=n^{-1}\sum_{i=1}^{n}{\mathbf{x_i}\mathbf{x_i^\prime}}$ $\mathbf{E}\psi^\prime(\epsilon_i)$ and $\mathbf{A_{2n}}=n^{-1}\sum_{i=1}^{n}{\mathbf{x_i}\mathbf{x_i^\prime}}\mathbf{E}\psi^2(\epsilon_i)$. \\

\hspace*{-6mm}(A.2)(iii)$^{\prime \prime}$ $n^{-2}\sum_{i = 1}^{n}||\mathbf{x_i}||^{12} + \sum_{i=1}^{n}\big[||\mathbf{\tilde{z}_i}||^4\max \{1, \mathbf{E}|\psi^\prime(\epsilon_i)|^4\}\big]= O(n^{-1})$.

\hspace*{-6mm}(A.3)(i)$^{\prime \prime}$ $\hspace*{1mm}$ $\mathbf{E}\psi(\epsilon_i) = 0$ for all $i\in \{1,\dots,n\}$.

\hspace*{-6mm}(A.3)(ii)$^{\prime \prime}$$\hspace*{0.5mm}$ $n^{-1}\sum_{i=1}^{n}\big[\mathbf{E}|\psi(\epsilon_i)|^{6+\upsilon}+\mathbf{E}|\psi^\prime(\epsilon_i)|^{6+\upsilon}+\mathbf{E}|\psi^{\prime\prime}(\epsilon_i)|^{4+\upsilon}\big]=O(1)$ for some $\upsilon>0$.

\hspace*{-6mm}(A.6)(i)$^{\prime \prime}$ $\hspace*{1mm}$ $\big(\psi(\epsilon_n), \psi'(\epsilon_n), \psi^2(\epsilon_n)\big)_{n=1}^{\infty}$ satisfies Cramer's condition in a uniform sense i.e. for any positive b,
\begin{align*}
\limsup_{n\rightarrow \infty}\sup_{||(t_1,t_2,t_3)||>b}\Big|\mathbf{E}\big(exp\big(it_1\psi(\epsilon_n)+it_2\psi'(\epsilon_n)+it_3\psi^2(\epsilon_n)\big)\big)\Big|<1.
\end{align*}

\hspace*{-6mm}(A.8) \hspace*{8mm}$\mathbf{A_{1n}}$ and $\mathbf{A_{2n}}$ both converge to non-singular matrices as $n\rightarrow \infty$.\\

We will denote the assumptions (A.1)-(A.4) by (A.1)$''$-(A.4)$''$ when (A.2) is defined with (iii)$''$ instead of (iii) and (A.3) is defined with (i)$^{\prime\prime}$, (ii)$^{\prime\prime}$ in place of (i) and (ii) respectively.

\subsection{Rate of Perturbation Bootstrap Approximation}
Note that when the regression errors are non-identically distributed, $\sqrt{n}\bm{\bar{\beta}_n}$ has asymptotic variance $\mathbf{A_{1n}^{-1}}$ $\mathbf{A_{2n}}\mathbf{A_{1n}^{-1}}$. Hence, the natural way of defining studentized pivot corresponding to $\bm{\bar{\beta}_n}$ is 
\begin{equation*}
\mathbf{\breve{H}_n}=\sqrt{n}\mathbf{\bar{\Sigma}_n^{-1/2}}(\bm{\bar{\beta}_n}-\bm{\beta})
\end{equation*}
where $\mathbf{\bar{\Sigma}_n^{-1/2}}=\mathbf{\bar{A}_{2n}^{-1/2}}\mathbf{\bar{A}_{1n}}$ with $\mathbf{\bar{A}_{1n}}=n^{-1}\sum_{i=1}^{n}\mathbf{x_ix'_i}\psi'(\bar{\epsilon}_i)$, $\mathbf{\bar{A}_{2n}}=n^{-1}\sum_{i=1}^{n}$ $\mathbf{x_ix'_i}\psi^2(\bar{\epsilon}_i)$ and $\bar{\epsilon}_i=y_i-\mathbf{x^\prime_i}\bm{\bar{\beta}_n}$, $i\in\{1,\dots,n\}$. Define the corresponding bootstrap pivot as  
\begin{equation*}
\mathbf{\breve{H}_n^*}=\sqrt{n}\mathbf{\Sigma_n^{*-1/2}}(\bm{\beta_n^*}-\bm{\bar{\beta}_n})
\end{equation*}
where $\mathbf{\Sigma_n^{*-1/2}}=\mathbf{A_{2n}^{*-1/2}}\mathbf{A_{1n}^*}$ with $\epsilon^*_i=y_i-\mathbf{x^\prime_i}\bm{\beta^*_n}$, $\mathbf{A_{1n}^*}=n^{-1}\sum_{i=1}^{n}\mathbf{x_ix'_i}\psi'(\epsilon^*_i)G_i^*$ and $\mathbf{A_{2n}^*}=n^{-1}\sum_{i=1}^{n}\mathbf{x_ix'_i}$ $\psi^2(\epsilon^*_i)(G_i-\mu_{G^*})^2$, $i\in\{1,\dots,n\}$.

\begin{thm}
Suppose, the assumptions \emph{(A.1)}$^{\prime\prime}$\emph{-(A.4)}$^{\prime\prime}$ and \emph{(A.5)(i)} hold.
 \begin{itemize}
 \item[\emph{(a)}]  Then there exist constant $C_5>0$ and a sequence of Borel sets $\mathbf{Q_{5n}}\subseteq \mathscr{R}^n$, such that $\mathbf{P}((\epsilon_1,\ldots,\epsilon_n)\in \mathbf{Q_{5n}})\rightarrow 1$ as $n\rightarrow\infty$, and given $(\epsilon_1,\ldots,\epsilon_n)\in \mathbf{Q_{5n}}$, $n\geq C_5$ such that there exists a sequence of statistics $\{\bm{\beta_n^*}\}_{n\geq 1}$ such that
\begin{equation*}
\mathbf{P_*}\big(\bm{\beta_n^*} \; solves  \;(2.1)\; and \;||\bm{\beta_n^*} - \bm{\bar{\beta}_n}||\leq C_5.n^{-1/2}.(log n)^{1/2}\big) \geq 1 - o\big(n^{-1/2}\big)
\end{equation*} 
 
\item[\emph{(b)}]  Suppose in addition \emph{(A.5)(ii),}\emph{(iii)$^\prime$} and $\emph{(A.8)}$ hold. Then there exist polynomial $\breve{a}^*_n(\cdot,\psi,G^{*})$ depending on first three moments of $G_1^{*}$ and on $\psi(\cdot)$, $\psi^\prime(\cdot)$ \& $\psi^{\prime\prime}(\cdot)$ through the residuals $\{\bar{\epsilon}_1,\dots,\bar{\epsilon}_n\}$, such that given $(\epsilon_1,.....,\epsilon_n)\in \mathbf{Q_{5n}}$, we have for $n\geq C_5$,
\begin{equation*}
\sup\limits_{B \in \mathscr{B}} |\mathbf{P_*}(\mathbf{\breve{H}_n^*}\in B) - \int_B\breve{\xi}^*_n(\mathbf{x})d\mathbf{x}| \leq \delta_n n^{-1/2}
\end{equation*}
where $\breve{\xi}^*_n(\mathbf{x}) = (1 + n^{-1/2}\breve{a}^*_n(\mathbf{x},\psi,G^{*}))\phi(\mathbf{x})$ and $\delta_n \equiv \delta_n(\epsilon_1,\ldots,\epsilon_n)$ tends to 0.\\

\item[\emph{(c)}]
Suppose, in addition to the assumptions \emph{(A.1)}$^{\prime\prime}$\emph{-(A.4)}$^{\prime\prime}$, \emph{(A.5)(i),(ii),}\emph{(iii)$^\prime$} and $\emph{(A.8)}$, \emph{(A.6)(i)$^{\prime\prime}$} holds. Then, for the collection of Borel sets defined by \emph{(3.1)},
\begin{equation*}
\sup\limits_{B \in  \mathscr{B}} \big{|}\mathbf{P_*}(\mathbf{\breve{H}_n^*}\in B) - \mathbf{P}(\mathbf{\breve{H}_n}\in B)\big{|} = o_p(n^{-1/2})
\end{equation*}
\end{itemize}

\end{thm}

\begin{rem}
The form of the studentized pivot $\mathbf{\breve{H}_n^*}$, defined for achieving second order correctness in non-IID case is different from $\mathbf{\tilde{H}_n^*}$, due to the difference in asymptotic variances of $\bm{\bar{\beta}_n}$ in two setups. In non-IID case, one cannot ignore computation of the negative square root of a matrix at each bootstrap iteration. But Theorem 5.1 is more general than Theorem 4.3 in the sense that it also includes the case when errors are IID. Note that $\mathbf{\bar{\Sigma}_n^*}=\mathbf{\bar{A}_{1n}^{*-1}}\mathbf{\bar{A}_{2n}^*}\mathbf{\bar{A}_{1n}^{*-1}}$ where $\mathbf{\bar{A}_{1n}^*}=n^{-1}\sum_{i=1}^{n}\mathbf{x_ix'_i}\psi'(\epsilon^*_i)$ and $\mathbf{\bar{A}_{2n}^*}=n^{-1}\sum_{i=1}^{n}\mathbf{x_ix'_i}\psi^2(\epsilon^*_i)$ and $\sigma_n^* = s_n^*\tau_n^{*-1}$ where $\tau_n^*=n^{-1}\sum_{i=1}^{n}\psi'(\epsilon_i^*)$, $s^{*2}_n = n^{-1}\sum_{i=1}^{n}\psi^2(\epsilon_i^*)$. We need to modify $\mathbf{\bar{\Sigma}_n^*}$ and ${\sigma}_n^*$ to $\mathbf{\Sigma_n^*}$ and $\tilde{\sigma}_n^*$ respectively to achieve second order correctness.
\end{rem}

\begin{rem}
There is no difference in employing perturbation bootstrap and the usual residual bootstrap with respect to the accuracy of inference. Under some mild conditions, both are second order correct. But in view Theorem 5.1, the advantage of employing perturbation bootstrap instead of residual counterpart is evident when the errors are no longer identically distributed. Perturbation bootstrap continues to be S.O.C. in non-IID case without any modification, whereas a modification in the resampling stage is required for residual bootstrap to achieve the same. To see this, consider the heteroscedastic simple linear regression model
\begin{align}\label{eqn:1.90}
y_i=\beta x_i + \epsilon_i
\end{align}
where $\epsilon_i$'s are independent, $\mathbf{E}\epsilon_i=0$ and $\mathbf{E}\epsilon_i^2=\sigma_i^2$. The least square estimator of $\beta$ is $\hat{\beta}=\sum_{i=1}^{n}x_iy_i/\sum_{i=1}^{n}x_i^2$ and hence $\mathbf{Var}(\hat{\beta})=\sum_{i=1}^{n}x_i^2\sigma_i^2/(\sum_{i=1}^{n}x_i^2)^2$. The bootstrap observations in residual bootstrap are $y_i^{**}=x_i \hat{\beta}+e_i^*$ where $\{e_1^*,\dots,e_n^*\}$ is a random sample from $\{(e_1-\bar{e}),\dots,(e_n-\bar{e})\}$, $\bar{e}=n^{-1}\sum_{i=1}^{n}e_i$ and $e_i=y_i-x_i \hat{\beta}$, $i\in \{1,\dots,n\}$, are least square residuals. The residual bootstrapped least square estimator is $\hat{\beta}^{**}=\sum_{i=1}^{n}x_iy_i^{**}/\sum_{i=1}^{n}x_i^2$. Hence, $\mathbf{Var}(\hat{\beta}^{**}|\epsilon_1,\dots,\epsilon_n)=\sum_{i=1}^{n}(e_i-\bar{e})^2/\sum_{i=1}^{n}x_i^2$ where $n^{-1}\sum_{i=1}^{n}[(e_i-\bar{e})^2-\sigma_i^2]\rightarrow 0$ as $n\rightarrow \infty$. Thus $\mathbf{Var}(\hat{\beta}^{**}|\epsilon_1,\dots,\epsilon_n)$ is not a consistent estimator of $\mathbf{Var}(\hat{\beta})$and hence residual bootstrap is not second order correct in approximating the distribution of $\hat{\beta}$ when errors are heteroscedastic. For details see Liu (1988). On the other hand, if $\hat{\beta}^*$ is the perturbation bootstrapped least square estimator, then it is easy to show 
$\mathbf{Var}(\hat{\beta}^{*}|\epsilon_1,\dots,\epsilon_n)=\sum_{i=1}^{n}x_i^2\sigma_i^2/(\sum_{i=1}^{n}x_i^2)^2 + O_p(n^{-1})$. Additionally, a centering adjustment is required in the definition of residual bootstrapped version of the regression M-estimator to achieve second order correctness even when the regression errors are IID [cf. Lahiri (1992)]; whereas in the perturbation bootstrap no adjustment is needed.
\end{rem}

\begin{rem}
In view of second order correctness of bootstrap in heteroscedastic linear regression, Theorem 5.1 is the most general result available. Nonparametric or residual bootstrap fails in heteroscedasticity, as shown by Liu (1988). Liu (1988) developed a weighted bootstrap method as a modification of residual bootstrap in least square setup for the simple linear regression model (5.1). She proposed the weight to be $x_i/\sum_{i=1}^{n}x_i^2$ corresponding to $i$th centered residual $(e_i-\bar{e}_n)$, $i\in \{1,\dots,n\}$, to achieve second order correctness. There is no general theory available on weighted bootstrap for the multiple linear regression model (1.1) even in heteroscedastic least square setup, to the best our knowledge. 
\end{rem}

\section{Proofs}

First we define some notations. Throughout this section, $C, C_1, C_2,\ldots$ will denote generic constants that do not depend on the variables like $n, x$, and so on. For a non-negative integral vector $\bm{\alpha} = (\alpha_1, \alpha_2,\ldots,\alpha_l)'$ and a function $f = (f_1,f_2,\ldots,f_l):\ \mathscr{R}^l\ \rightarrow \ \mathscr{R}^l$, $l\geq 1$, write $|\bm{\alpha}| = \alpha_1 +\ldots+ \alpha_l$, $\bm{\alpha}! = \alpha_1!\ldots \alpha_l!$, $f^{\bm{\alpha}} = (f_1^{\alpha_1})\ldots(f_l^{\alpha_l})$. For $\mathbf{t} =(t_1,\ldots t_l)'\in \mathscr{R}^l$ and $\mathbf{\alpha}$ as above, define $t^{\bm{\alpha}} = t_1^{\alpha_1}\ldots t_l^{\alpha_l}$. The collection $\mathscr{B}$ will always be used to denote the collection of Borel subsets of $\mathscr{R}^p$ which satisfy (3.1). $\mu_{G^*}$ and $\sigma_{G^*}^2$ will respectively denote mean and variance of $G_1^*$. We want to mention here that only the important steps are presented in the proofs of the proposition and the theorems. For further details see the supplementary material Das and Lahiri (2017). Although the proofs for second order results of perturbation bootstrap go through more or less same line as that for residual bootstrap in Lahiri (1992), the advantage in perturbation bootstrap is that the perturbing quantities are independent of the regression errors and hence it is much easier to obtain suitable stochastic approximation to the bootstrapped pivot and finally the EE than the same in case of residual bootstrap. On the negative side, in our proofs atleast we need Cramer's condition separately on regression errors and on the perturbing quantities [see assumptions (A.5) and (A.6)], whereas for residual bootstrap, one can derive a restricted Cramer's condition on resampled residuals from the Cramer's condition on regression errors to obtain second order correctness. Moreover, second order results can be established for residual bootstrap, after a modification, without any Cramer type condition in the case $p=1$ [cf. Karabulut and Lahiri (1997)]. We do not know yet if similar conclusion can be drawn in case of perturbation bootstrap.\\ 

Before coming to the proofs we state some lemmas:

\begin{lemma}
Let, $\{\mathbf{Y_{i}}=(Y_{i1}, Y_{i2})', 1\leq i \leq n\}$  be a collection of mean zero independent random vectors. Define, for some non random vectors $\mathbf{l_{1i}}$ and $\mathbf{l_{2i}}$ of dimensions $p_1$ and $p_2$ respectively with $\sum_{i=1}^{n}\mathbf{l_{ji}l'_{ji}}=\mathbf{I_{p_j}}$ and $\tilde{\gamma}_n = (\sum_{j=1}^{2}\sum_{i=1}^{n}||\mathbf{l_{ji}}||^4)^{1/2} =O(n^{-1/2})$,
\begin{equation*}
\mathbf{U_{i}} = (\mathbf{l'_{1i}}Y_{i1},\mathbf{l'_{2i}}Y_{i2})',\;\;\;\;\;  \mathbf{V_n} = \mathbf{Cov}\Big(\sum_{i = 1}^{n}\mathbf{U_{i}}\Big), \;\;\;\;\;\; \mathbf{\tilde{U}_{i}} = \mathbf{V_{n}^{-1/2}}\mathbf{U_{i}}
\end{equation*}
for $1\leq i\leq n$, and $\mathbf{S_n} = \sum_{i=1}^{n}\mathbf{\tilde{U}_{i}}$. Let $\mathbf{\tilde{\alpha}_n} = n^{-1}\sum_{i=1}^{n}\mathbf{E}||\mathbf{Y_{i}}||^3I(||\mathbf{Y_{i}}||^2>\lambda 
\tilde{\gamma}_n^{-1})$, where $I(\cdot)$ is the indicator function and $\lambda$ satisfies $0<\lambda< \liminf\limits_{n\rightarrow\infty}\lambda_n$, $\lambda_i =$ the smallest eigen value of $\mathbf{\Sigma_i}$, $\mathbf{\Sigma_i} = \mathbf{Cov}(\mathbf{Y_{i}})$.
Suppose, $\{\mathbf{M_{0n}}\}_{n\geq 1}$, $\{\mathbf{M_{in}}\}_{n\geq 1}, i = 1,\ldots,p$ be $(p+1)$ sequence of matrices such that for each $n\geq 1$, $\mathbf{M_{0n}}$ is of order $p\times (p+r)$. and $\mathbf{M_{in}}, 1\leq i\leq p$, are of order $(p+r)\times (p+r)$, $p\geq 1$, $r\geq 1$. Let, $k=p+r$, $\mathbf{\bar{M}_{0n}} = [\mathbf{0}: I_r]_{r\times k}$ and $\mathbf{\tilde{M}_{0n}}=[\mathbf{M_{0n}}':\mathbf{\bar{M}_{0n}}']'$. Define the functions $g_n:\mathscr{R}^k\rightarrow \mathscr{R}^p$ by $g_n(\mathbf{x})=\mathbf{M_{0n}}\mathbf{x} + (\mathbf{x}'\mathbf{M_{1n}}\mathbf{x},\ldots,\mathbf{x}'\mathbf{M_{pn}}\mathbf{x})'$, $\mathbf{x}\in \mathscr{R}^k$, $n\geq 1$. Assume that
\begin{enumerate}
\item[\emph{(a)}] there exists a constant $k$ such that $n^{-1}\sum_{i=1}^{n}\mathbf{E}||\mathbf{Y_{i}}||^3 < k$ for all $n\geq 1$.
\item[\emph{(b)}] $\tilde{\alpha}_n = o(1)$.
\item[\emph{(c)}] the characteristic function $g_n$ of $\mathbf{Y_{n}}$ satisfies $\limsup_{n\rightarrow \infty}\sup_{||(\bm{t})||>b}|g_n(\bm{t})|<1$ for all $b>0$.
\item[\emph{(d)}] $max\{||\mathbf{M_{in}}||: 1\leq i\leq p\} = O(\tilde{\gamma}_n)$.
\item[\emph{(e)}] $||\mathbf{M_{0n}}|| = O(1)$, $\liminf_{n\rightarrow\infty} \inf\{||\mathbf{\tilde{M}_{0n}u}||: ||\mathbf{u}|| = 1, \mathbf{u}\in \mathscr{R}^k\} \geq \delta$ for some constant $\delta>0$.
\end{enumerate}
Then for the class $\mathscr{B}$ of Borel sets satisfying (3.1),
\begin{equation*}
\sup\limits_{B \in \mathscr{B}} \big|\mathbf{P}(g_n(\mathbf{S_n})\in B) - \int_B\mathring{\xi}_n(\mathbf{x})d\mathbf{x}\big| = o(\tilde{\gamma}_n) \;\;\;\; as \; n\rightarrow \infty
\end{equation*}
where $\mathring{\xi}_n(.) = (1 + n^{-1/2}\mathring{a}(\cdot))\phi_{\mathbf{\mathring{D}_n}}(\cdot)$, $\mathbf{\mathring{D}_n} = \mathbf{M_{0n}}\mathbf{M_{0n}^\prime}$ and $\mathring{a}(\cdot)$ is a polynomial whose coefficients are continuous functions of $\mathbf{E}(\mathbf{Y_{i}})^{\bm{\alpha}}, |\mathbf{\bm{\alpha}}|\leq 3$ and $i\in\{1,\dots,n\}$.
\end{lemma}

\subsubsection*{proof :} 
The above Lemma follows from Theorem 20.6 of Bhattacharya and Rao (1986) and retracting the proofs of Lemma 3.1 and 3.2 of Lahiri (1992).

\begin{lemma}
Under the assumptions \emph{(A.1)}-\emph{(A.3)} or \emph{(A.1)$^{\prime\prime}$}-\emph{(A.3)$^{\prime\prime}$}, it follows that \\
$\big(\sum_{i=1}^{n}||\mathbf{\breve{v}_i}||^4\big)^{1/2}= O_p(n^{-1/2})$.
\end{lemma}
\subsubsection*{proof :}
See supplementary material Das and Lahiri (2017).

\begin{lemma}
Under the assumptions \emph{(A.2) (i)} and \emph{(A.2) (iii)} or \emph{(A.2) (iii)$^{\prime\prime}$}, the following is true.
\begin{enumerate}
\item[\emph{(a)}] $\big(\sum_{i=1}^{n}||\mathbf{d_i}||^6\big)^{1/4} + \big(\sum_{i=1}^{n}||\mathbf{d_i}||^4\big)^{1/2}= O(n^{-1/2})$.
\item[\emph{(b)}] $\sum_{i=1}^{n}||\mathbf{x_i}||^{j} = O(n)$ for $j =3,4,5,6,6+2\alpha$ when the errors are IID and for $j=6+2\alpha, 3,\dots, 12$ when the errors are non-IID.
\end{enumerate}
\end{lemma}
\subsubsection*{proof :} 
This lemma follows from assumption (A.2) and by applying H\"{o}lders inequality.

\vspace{5mm}

We present only outline of the proofs of the main results from Section 4 and 5 to save space. For details, see the supplementary material Das and Lahiri (2017).

\subsection{Outline of the proof of Proposition 4.1}
Suppose,
\begin{equation*}
\sum_{i = 1}^{n}\mathbf{x_i}\psi(y_i - \mathbf{x'_i}\mathbf{t_n^*})G^*_i = \mathbf{0}
\end{equation*}
Then by Taylor's expansion we have,
\begin{equation}
\sum_{i = 1}^{n}\mathbf{x_i}\psi(\bar{\epsilon}_i)G_i^* + \sum_{i = 1}^{n}\mathbf{x_ix_i'}(\bm{\bar{\beta}_n} - \mathbf{t_n^*})\psi'(\bar{\epsilon}_i)G_i^* + \sum_{i = 1}^{n} \mathbf{x_i}\dfrac{[\mathbf{x_i}'(\bm{\bar{\beta}_n} - \bm{t_n^*})]^2}{2}\psi''(u_i)G_i^* = 0
\end{equation}
where for each $i\in \{1,\ldots,n\}$, $|u_i - \bar{\epsilon}_i|\leq |\epsilon_i^* - \bar{\epsilon}_i|$.\\
Now (6.1) can be written as
\begin{equation}
L_n^*(\mathbf{t_n^*} - \bm{\bar{\beta}_n}) = \Delta_n^* + R_n^*
\end{equation}
where\\ 
$\Delta_n^* = n^{-1}\sum_{i = 1}^{n}\mathbf{x}_i\psi(\bar{\epsilon}_i)(G_i^* - \mu_{G^*})$\\
$L_n^* = n^{-1}\sum_{i = 1}^{n}\mathbf{x}_i\mathbf{x}_i'\psi'(\bar{\epsilon}_i)G_i^*$\\
$\mathbf{E_*}L_n^* = n^{-1}\sum_{i = 1}^{n}\mathbf{x}_i\mathbf{x}_i'\psi'(\bar{\epsilon}_i)\mu_{G^*}$\\
$R_n^* = n^{-1}\sum_{i = 1}^{n}\mathbf{x}_i\dfrac{[\mathbf{x}_i'(\bm{\bar{\beta}}_n - \mathbf{t}_n^*)]^2}{2}\psi''(u_i)G_i^*$
\vspace*{2mm}

By Fuk and Nagaev inequality (1971) [hereafter referred to as FN(71)], lemma 6.3, the Lipschitz property of $\psi''(\cdot)$ and the Taylor's expansion of $\psi(\cdot)$ and $\psi'(\cdot)$, it follows that
there exist a constant $C>0$ and a sequence of Borel sets $\mathbf{Q}_{n}\subseteq \mathscr{R}^n$, such that given $(\epsilon_1,.....,\epsilon_n)\in \mathbf{Q}_{n}$ with $\mathbf{P}((\epsilon_1,......,\epsilon_n)\in \mathbf{Q}_{n})\rightarrow 1$ , for $n \geq C$
and any $0<\epsilon<1$, 
\begin{equation}
\mathbf{P_*}\Big( \big{|}\sum_{i=1}^{n}||\mathbf{x_i}||^{3 + \alpha}(G_i^* - EG_i^*)\big{|}>n\epsilon\Big) = o(n^{-1/2})
\end{equation}

\begin{align}
\mathbf{P_*}\Big(\big{|}\sum_{i = 1}^{n}x_{ij}x_{ik}\psi'(\bar{\epsilon}_i)(G_i^* - EG_i^*)\big{|}>n{\epsilon}\Big)=o(n^{-1/2}),\;\;\;\; j,k\in\{1,\ldots,p\}
\end{align}

\begin{equation}
\mathbf{P_*}\Big(||\Delta_n^*||>C.n^{-1/2}(log n)^{1/2}\Big) = o(n^{-1/2})
\end{equation}
\hspace{2mm}

Hence, from (6.3)-(6.5), on the set $\mathbf{Q_n}$ and given $(\epsilon_1,.....,\epsilon_n)\in \mathbf{Q_{n}}$ with $\mathbf{P}((\epsilon_1,......,\epsilon_n)$ $\in \mathbf{Q_{n}})\rightarrow 1$, for $n\geq C_1$, (6.2) can be rewritten as 
$(\mathbf{t_n^*} - \bm{\bar{\beta}}_n) = f_n(\mathbf{t}_n^* - \bm{\bar{\beta}}_n)$, 
where $f_n$ is a continuous function from $\mathscr{R}^p$ to $\mathscr{R}^p$ satisfying $\mathbf{P_*}(||f_n(\mathbf{t_n}^* - \bm{\bar{\beta}}_n)||\leq C_1.n^{-1/2}(logn)^{1/2}) = 1 - o(n^{-1/2})\;\;\; as\; n\rightarrow\infty$ whenever $||\mathbf{t_n^*} - \bm{\bar{\beta}_n}||\leq C_1.n^{-1/2}(logn)^{1/2}$ for some constants $C_1>0$.

Hence, Proposition 4.1 follows by Brouwer's fixed point theorem.

\subsection{Outline of the proof of Theorem 4.1}

Consider, the sequence of statistics $\{\bm{\beta_n^*}\}_{n\geq 1}$ which satisfies the proposition. Then (6.2) can be written as
\begin{align}
\sqrt{n}(\bm{\beta_n^*} - \bm{\bar{\beta}_n}) &= L_n^{*-1}\sqrt{n}[\Delta_n^* + \tilde{\chi}_n^* + R_{1n}^*]\\
& = L_n^{*-1}\sqrt{n}\Delta_n^* + R_{2n}^*
\end{align}
where $\;\;\tilde{\chi}_n^* = n^{-1}\sum_{i = 1}^{n}\mathbf{x_i}\dfrac{[\mathbf{x_i}'(\bm{\beta_n^* }- \bm{\bar{\beta}_n})]^2}{2}\psi''(\bar{\epsilon}_i)G_i^*$\\ 

Now, by FN(71), for some constant $C>0$,

\begin{equation*}
\mathbf{P_*}(||R_{1n}^*||>C.n^{-(2+\alpha)/2}(logn)^{(2+\alpha)/2}) = o_p(n^{-1/2})
\end{equation*}
and
\begin{equation*} 
\mathbf{P_*}(||R_{2n}^*||>C.n^{-1/2}(log n)) = o_p(n^{-1/2})
\end{equation*}

Again,
\begin{equation}
L_n^{*-1} = (\mathbf{E_*}L_n^*)^{-1} + W_n^* + \tilde{Z}_n^*
\end{equation}
where\\

$W_n^* = (\mathbf{E_*}L_n^*)^{-1}(\mathbf{E_*}L_n^* - L_n^*)(\mathbf{E_*}L_n^*)^{-1}$ 

$\tilde{Z}_n^* = (\mathbf{E_*}L_n^*)^{-1}(\mathbf{E_*}L_n^* -L_n^*)(\mathbf{E_*}L_n^*)^{-1}(\mathbf{E_*}L_n^* -L_n^*)L_n^{*-1}
$\\

Now, it can be shown by FN(71) that for some constant $C_1>0$, as $n\geq C_1$,
\begin{align}
\mathbf{P_*}(||\tilde{Z}_n^*||&>C_1.n^{-1/2}(logn)^{-1})\nonumber \\
&\leq \mathbf{P_*}(||L_n^* - \mathbf{E_*}L_n^*||>C_1.n^{-1/4}(log n)^{-1/2})\nonumber \\
&= o_p(n^{-1/2})
\end{align}

Therefore, it follows that there exists $C_2>0$ and a sequence of Borel sets $\mathbf{Q_{2n}}$, such that $\mathbf{P}((\epsilon_1,......,\epsilon_n)\in \mathbf{Q_{2n}})\rightarrow 1$ as $n\rightarrow\infty$, and given $(\epsilon_1,.....,\epsilon_n)\in \mathbf{Q_{2n}}$ and $n\geq C_2$,

\begin{equation}
\sqrt{n}(\bm{\beta_n^*} -\bm{\bar{\beta}_n})= (\mathbf{E_*}L_n^*)^{-1}\sqrt{n}\Delta_n^* +  W_n^*\sqrt{n}\Delta_n^*+  (\mathbf{E_*}L_n^*)^{-1}\sqrt{n}\chi_n^{*} + R_{3n}^*
\end{equation}
where 
$\;\chi_n^{*} = n^{-1}\sum_{i = 1}^{n}\mathbf{x_i}\dfrac{[\mathbf{x_i}'((\mathbf{E_*}L_n^*)^{-1}\Delta_n^*)]^2}{2}\psi''(\bar{\epsilon}_i)\mu_{G^*}$\\
and 
\begin{equation*}
\mathbf{P_*}(||R_{3n}^*||=o(n^{-1/2})) = 1-o(n^{-1/2})
\end{equation*}

Since $\mathbf{\bar{\Sigma}_n^{-1/2}}=O_p(1)$, so by argument similar to (4.12) of Qumsiyeh (1990a), we have
\begin{align}
\sup\limits_{B \in  \mathscr{B}} |\mathbf{P_*}(\mathbf{F_n^*}\in B ) - \mathbf{P_*}(\mathbf{U_n^*} \in B )| = o_p(n^{-1/2})
\end{align}
where $\mathbf{U_n^*} = \sqrt{n}\mathbf{\bar{\Sigma}_n^{-1/2}}\bigg[(\mathbf{E_*}L_n^*)^{-1}\Delta_n^* +  W_n^*\Delta_n^*+  (\mathbf{E_*}L_n^*)^{-1}\chi_n^{*}\bigg]$\\

Now, for all $1\leq i\leq n$, defining
$Y_{i}^* = (G_i^* - \mu_{G^*})$,
$\mathbf{X_i^*} = \mathbf{\breve{v}_i}Y^*_{i}$,
$\mathbf{V_n^*} = \sum_{i=1}^{n}\mathbf{Cov_*}(\mathbf{X_i^*})$,
$\mathbf{\tilde{X}_i^*} = \mathbf{V_n^{*-1/2}X_i^*}$ and
$\mathbf{S_n^*} = \sum_{i = 1}^{n}\mathbf{\tilde{X}_i^*}$,
it can be established that

\begin{equation}
\mathbf{U_n^{*}} = \mathbf{M_{0n}^*S_n^*} + (\mathbf{{S^*_n}}' \mathbf{M_{1n}^*S_n^*},\ldots, \mathbf{{S^*_n}}' \mathbf{M_{pn}^*S_n^*})'
\end{equation}
where $\mathbf{M_{0n}^*} = O_p(1)$ and $\mathbf{M_{jn}^*} = O_p(n^{-1/2})$ for all $j\in \{1,\ldots,p\}$.\\

Therefore, by Lemma 6.1 and 6.2,
\begin{equation}
\sup\limits_{B \in  \mathscr{B}} |\mathbf{P_*}(\mathbf{U_n^*}\in B) - \int_B\xi^*_n(\mathbf{x})d\mathbf{x}| = o_p(n^{-1/2}) \;\;\;\; as \; n\rightarrow \infty
\end{equation}
where 
\begin{equation}
{\xi}^*_n(\mathbf{x}) = \Bigg{[}1-n^{-1/2}\Big{\{}{\sum_{|\bm{\nu}|=1}b_{11}^{*(\bm{\nu})}D^{\bm{\nu}}}+\sum_{|\bm{\nu}|=3}\dfrac{b_{31}^{*(\bm{\nu})}}{\bm{\nu}!}D^{\bm{\nu}}\Big{\}}\Bigg{]}\phi(\mathbf{x})
\end{equation}

Now, the coefficients $b_{11}^{*(\bm{\nu})}$ and $b_{31}^{*(\bm{\nu})}$ can be computed using the transformation techniques of Bhattacharya and Ghosh (1978). If $\bm{\nu_1}$ is a $p\times 1$ vector with all the elements being 0, except the $j$th one and $\bm{\nu_2}$ is a $p\times 1$ vector with all the elements being 0, except the $j_1,j_2$ and $j_3$ positions then after some algebraic calculations it can be shown that

\begin{align}
b_{11}^{*(\bm{\nu}_1)} = &\sum_{k=1}^{p}h_{jkn}\Big(n^{-1}\sum_{i=1}^{n}\big[\mathbf{z'_iE^*_{kn}\bar{A}_{1n}^{-1}x_i}\psi(\bar{\epsilon}_i)\psi'(\bar{\epsilon}_i)\big]\Big) \nonumber\\
&+ (2n)^{-1}\sum_{i=1}^{n}a_{jin}^*\mathbf{x'_i\bar{A}_{1n}^{-1}\bar{A}_{2n}\bar{A}_{1n}^{-1}x_i}\psi''(\bar{\epsilon}_i)
\end{align}
\begin{align}
b_{31}^{*(\bm{\nu}_2)} = 
 & n^{-1}\sum_{i=1}^{n}\Bigg[\bigg(\prod_{m=1}^{3}a_{j_{m}in}^*\bigg)\psi^3(\bar{\epsilon}_i)\Big] \nonumber \\
& +2n^{-2}\sum_{i,j=1}^{n}\Big[a_{j_1in}^*a_{j_2in}^*\big(\sum_{k=1}^{p}h_{j_3kn}\mathbf{z'_iE_{kn}^*\bar{A}_{1n}^{-1}x_j}\big)\psi^2(\bar{\epsilon}_i)\psi(\bar{\epsilon}_j)\psi'(\bar{\epsilon}_j)\Big] \nonumber\\
&+2n^{-2}\sum_{i,j=1}^{n}\Big[a_{j_1in}^*a_{j_3in}^*\big(\sum_{k=1}^{p}h_{j_2kn}\mathbf{z'_iE_{kn}^*\bar{A}_{1n}^{-1}x_j}\big)\psi^2(\bar{\epsilon}_i)\psi(\bar{\epsilon}_j)\psi'(\bar{\epsilon}_j)\Big] \nonumber
\end{align}
\begin{align}
\hspace*{18mm}&+2n^{-2}\sum_{i,j=1}^{n}\Big[a_{j_2in}^*a_{j_3in}^*\big(\sum_{k=1}^{p}h_{j_1kn}\mathbf{z'_iE_{kn}^*\bar{A}_{1n}^{-1}x_j}\big)\psi^2(\bar{\epsilon}_i)\psi(\bar{\epsilon}_j)\psi'(\bar{\epsilon}_j)\Big] \nonumber \\
&+3n^{-3}\sum_{i,j,l=1}^{n}a_{j_1in}^*a_{j_2in}^*a_{j_3in}^*\big(\mathbf{x'_j\bar{A}_{1n}^{-1}x_lx'_l\bar{A}_{1n}^{-1}x_i}\big)\psi''(\bar{\epsilon}_l)\psi^2(\bar{\epsilon}_i)\psi^2(\bar{\epsilon}_j)
\end{align}

where $\mathbf{\bar{A}_{1n}}$ and $\mathbf{\bar{A}_{2n}}$ are as defined earlier and
$\mathbf{\bar{A}_{2n}^{-1/2}}=(\mathbf{h_{1n}},\ldots,\mathbf{h_{pn}})$, $\;\;\mathbf{h'_{jn}x_i}=a_{jin}^*$,  $\mathbf{h_{jn}}=(h_{1jn},\ldots$ $,h_{pjn})$, $j\in\{1,\cdots,p\}, i\in\{1,\ldots,n\}$ and $\mathbf{E_{kn}^*}$ is a $q \times p$ matrix with $||\mathbf{E_{kn}^*}||\leq q$ for all $k\in \{1,\ldots,p\}$\\

Now, one can find the two term EE of $\mathbf{F_n}=\sqrt{n}\sigma^{-1}\mathbf{A_n^{1/2}}(\bm{\bar{\beta}_n}-\bm{\beta})$ in similar way such that (for detail see Lahiri(1992)) 
\begin{equation}
\sup\limits_{B \in  \mathscr{B}} |\mathbf{P}(\mathbf{F_n}\in B) - \int_B\xi_n(\mathbf{x})d\mathbf{x}| = o(n^{-1/2}) \;\;\;\; as \; n\rightarrow \infty
\end{equation}
where 
\begin{equation}
\xi_n(\mathbf{x}) = \Bigg{[}1-n^{-1/2}\Big{\{}{\sum_{|\bm{\nu}|=1}b_{11}^{(\bm{\nu})}D^{\bm{\nu}}}+\sum_{|\bm{\nu}|=3}\dfrac{b_{31}^{(\bm{\nu})}}{\bm{\nu}!}D^{\bm{\nu}}\Big{\}}\Bigg{]}\phi(\mathbf{x})
\end{equation}
where the coefficients $b_{11}^{(\bm{\nu_1})}$ and $b_{31}^{(\bm{\nu_2})}$ are such that for all $j,j_1,j_2,j_3 \in \{1,\ldots,p\}$,
$\big(b_{11}^{*(\bm{\nu_1})}-b_{11}^{(\bm{\nu_1})}\big)$ and $\big(b_{31}^{*(\bm{\nu_2})}-b_{31}^{(\bm{\nu_2})}\big)$ both can be shown to converge in probability to $0$. Hence by (6.12)-(6.18), Theorem 4.1 follows.

\subsection{Outline of the proof of Theorem 4.2}

\begin{flushleft}
We have,
\end{flushleft}
\begin{equation}
\mathbf{H^*_n} = \sqrt{n}\sigma_n^{*-1}\hat{\sigma}_n\mathbf{\bar{\Sigma}_n^{-1/2}}(\bm{\beta^*_n} - \bm{\bar{\beta}_n})
\end{equation}
where $\sigma_n^*$ is as defined earlier.
Now using Taylor's expansion and Lipschitz property of $\psi''(\cdot)$, it can be established that 

\begin{equation}
\mathbf{H^*_n} = \mathbf{F_n^*} - \sqrt{n}\hat{\sigma}_n\mathbf{\bar{\Sigma}_n^{-1/2}}Z^*_n((\mathbf{E_*}L_n^*)^{-1}\Delta_n^*) +R^*_{4n}
\end{equation}
where 
\begin{align*}
Z_n^*=&(2s^3_n|\tau_n|)^{-1}\Big[2\tau_n s_n^2 \Big(\dfrac{1}{n}\sum_{i=1}^{n}\psi''(\bar{\epsilon}_i) [\mathbf{x_i}'((\mathbf{E_*}L_n^*)^{-1}\Delta_n^*)]\Big)\\ 
&- \tau^2_n\Big(\dfrac{2}{n}\sum_{i=1}^{n}\psi(\bar{\epsilon}_i) \psi'(\bar{\epsilon}_i)[\mathbf{x_i}'((\mathbf{E_*}L_n^*)^{-1}\Delta_n^*)]\Big)\Big]
\end{align*}
and there exist constants $C_3>0$ and a sequence of Borel sets $\mathbf{Q_{3n}}$ such that $\mathbf{P}(\mathbf{Q_{3n}})\uparrow 1$ and given $(\epsilon_1,........\epsilon_n)\in \mathbf{Q_{3n}}$ and $n\geq C_3$,
\begin{equation}
\mathbf{P_*}(||R^*_{4n}||=o(n^{-1/2})) = 1-o(n^{-1/2})
\end{equation}

Therefore, writing $\mathbf{H_n^*}$ as $\mathbf{H_n^*} = \mathbf{\tilde{U}_n^*} + R_{4n}^*$, we have
\begin{equation}
\mathbf{\tilde{U}_n^{*}} = \mathbf{\tilde{M}_{0n}^*S_n^*} + (\mathbf{{S^*_n}}' \mathbf{\tilde{M}_{1n}^*S_n^*},\ldots, \mathbf{{S^*_n}}'\mathbf{\tilde{M}_{pn}^*S_n^*})'
\end{equation}
where $\mathbf{\tilde{M}_{0n}^*} = O_p(1)$ and $\mathbf{\tilde{M}_{jn}^*} = O_p(n^{-1/2})$ for all $j\in \{1,\ldots,p\}$.\\

Hence, by Lemma 6.1,
\begin{equation}
\sup\limits_{B \in  \mathscr{B}} |\mathbf{P_*}(\mathbf{\tilde{U}_n^*}\in B) - \int_B\tilde{\xi}^*_n(\mathbf{x})d\mathbf{x}| = o_p(n^{-1/2}) \;\;\;\; \text{as} \; n\rightarrow \infty
\end{equation}
where 
\begin{equation}
\tilde{\xi}^*_n(\mathbf{x}) = \Bigg{[}1-n^{-1/2}\Big{\{}{\sum_{|\bm{\nu}|=1}\tilde{b}_{11}^{*(\bm{\nu})}D^{\bm{\nu}}}+\sum_{|\bm{\nu}|=3}\dfrac{\tilde{b}_{31}^{*(\bm{\nu})}}{\bm{\nu}!}D^{\bm{\nu}}\Big{\}}\Bigg{]}\phi(\mathbf{x})
\end{equation}

Hence part (a) follows by (4.12) of Qumsiyeh (1990a).\\

Suppose the two term EE of the original studentized regression M-estimator $\mathbf{H_n}=\sqrt{n}\hat{\sigma}_n^{-1}\mathbf{A_n^{1/2}}(\bm{\bar{\beta}_n}$ $-\bm{\beta})$ is  \\
\begin{equation}
\tilde{\xi}_n(\mathbf{x}) = \Bigg{[}1-n^{-1/2}\Big{\{}{\sum_{|\bm{\nu}|=1}\tilde{b}_{11}^{(\bm{\nu})}D^{\bm{\nu}}}+\sum_{|\bm{\nu}|=3}\dfrac{\tilde{b}_{31}^{(\bm{\nu})}}{\nu!}D^{\bm{\nu}}\Big{\}}\Bigg{]}\phi(\mathbf{x})
\end{equation}

Now part (b) of Theorem 4.2 follows directly by comparing (6.24) and (6.25). Again after some algebraic calculations, it can be shown that $\tilde{b}_{11}^{(\bm{\nu})}$ and $\tilde{b}_{31}^{(\bm{\nu})}$ both contain terms involving $\big[2E\psi^2(\epsilon_1)$ $E\psi(\epsilon_1)\psi'(\epsilon_1)- E\psi'(\epsilon_1)E\psi^3(\epsilon_1)\big]$ which cannot be replicated by the terms present in $\tilde{b}_{11}^{*(\bm{\nu})}$ and $\tilde{b}_{31}^{*(\bm{\nu})}$ [cf. Supplementary material Das and Lahiri (2017)]. Hence part (c) of Theorem 4.2 follows.

\subsection{Outline of the proof of Theorem 4.3}

\begin{flushleft}
We have the modified studentized bootstrapped M-estimator as,
\end{flushleft}
\begin{equation}
\mathbf{\tilde{H}^*_n} = \sqrt{n}(\tilde{\sigma}_n^*)^{-1}\hat{\sigma}_n\mathbf{\bar{\Sigma}_n^{-1/2}}(\bm{\beta^*_n} - \bm{\bar{\beta}_n})
\end{equation}
where $\tilde{\sigma}_n^*=\tilde{s}_n^*\tilde{\tau}_n^{*-1}$,  $\tilde{\tau}_n^*=n^{-1}\sum_{i=1}^{n}\psi'(\epsilon_i^*)G_i^*$ and $\tilde{s}_n^{*2}=n^{-1}\sum_{i=1}^{n}\psi^2(\epsilon_i^*)(G_i^*-\mu_{G^*})^2$. Also suppose, $\bar{\tau}_n = \mu_{G^*}\tau_n$ and $\bar{s}_n^2=\sigma_{G^*}^2s_n^2$.\\

Now using the same line of arguments which is working behind (6.20) in the proof of Theorem 4.2, it can be shown that 

\begin{equation}
\mathbf{\tilde{H}^*_n} = \mathbf{F_n^*} - \sqrt{n}\hat{\sigma}_n\mathbf{\bar{\Sigma}_n^{-1/2}}(Z_{n}^*-\bar{Z}_{n}^*)((\mathbf{E_*}L_n^*)^{-1}\Delta_n^*) +R^*_{5n}
\end{equation}
where $\bar{Z}_n^*$ is as defined in the proof of Theorem 4.2 and $\tilde{Z}_n^*$ is defined as
\begin{align*}
\bar{Z}_n^*=&2^{-1}\big(\bar{\tau}_n\bar{s}_n\big)^{-2}\Bigg[2\bar{\tau}_n \bar{s}_n^2\Big(n^{-1}\sum_{i=1}^{n}\psi'(\bar{\epsilon}_i)(G_i^*-\mu_{G^*})\Big)\\ 
&- \bar{\tau}_n^2 \Big(n^{-1}\sum_{i=1}^{n}\psi^2(\bar{\epsilon}_i)[(G_i^*-\mu_{G^*})^2-\sigma_{G^*}^2]\Big)\Bigg]
\end{align*}
and there exist constant $C_4>0$ and a sequence of Borel sets $\mathbf{Q_{4n}}$ such that $\mathbf{P}(\mathbf{Q_{4n}})\uparrow 1$ and given $(\epsilon_1,........\epsilon_n)\in \mathbf{Q_{4n}}$ and $n\geq C_4$,
\begin{equation}
\mathbf{P_*}(||R^*_{5n}||=o(n^{-1/2})) = 1-o(n^{-1/2})
\end{equation}

Therefore, defining $Y_{1i}^* = G_i^*-\mu_{G^*}$,
$Y_{2i}^* = (G_i^*-\mu_{G^*})^2-\sigma_{G^*}^2$\\
$\mathbf{X_i^*} = \Big(\mathbf{\breve{v}'_i}Y_{1i}^*, n^{-1/2}\psi^2(\bar{\epsilon}_i)Y_{2i}^*\Big)'$,
$\mathbf{V_n^*} = \sum_{i=1}^{n}\mathbf{Cov_*}(\mathbf{X_i^*})$,
$\mathbf{\tilde{X}_i^*} = \mathbf{V^{*-1/2}_nX_i^*}$,
$\mathbf{\bar{S}_n^*} = \sum_{i=1}^{n}\mathbf{\tilde{X}_i^*}$ with $\mathbf{\bar{v}_i}$ defined with $\mathbf{\breve{z}_i}$ in place of $\mathbf{z_i}$.\\

Hence, we have $\mathbf{\tilde{H}_n^*}$ as $\mathbf{\tilde{H}_n^*} = \mathbf{\bar{U}_n^*} + R_{5n}^*$, where
\begin{equation}
\mathbf{\bar{U}_n^{*}} = \mathbf{\bar{M}_{0n}^*}\mathbf{\bar{S}_n^*} + (\mathbf{{\bar{S}^{*'}_n}} \mathbf{\bar{M}_{1n}^*} \mathbf{\bar{S}_n^*},\ldots, \mathbf{{\bar{S}^{*'}_n}} \mathbf{\bar{M}_{pn}^*\bar{S}_n^*})'
\end{equation}
with $\mathbf{\bar{M}_{0n}^*} = O_p(1)$ and $\mathbf{\bar{M}_{jn}^*} = O_p(n^{-1/2})$ for all $j\in \{1,\ldots,p\}$.\\

Hence, there exists a two term EE $\bar{\xi}^*(\cdot)$, as in Theorem 4.2, such that
\begin{equation}
\sup\limits_{B \in  \mathscr{B}} |\mathbf{P_*}(\mathbf{\tilde{H}_n^*}\in B) - \int_B\bar{\xi}^*_n(\mathbf{x})d\mathbf{x}| = o_p(n^{-1/2}) \;\;\;\; as \; n\rightarrow \infty
\end{equation}
 
Now, $\bar{\xi}_n^*(\cdot)$ can be found explicitly as in standardized case. See supplementary material Das and Lahiri (2017) for more details. Again if $\bar{\xi}_n^*(\cdot)$ is compared with $\tilde{\xi}_n(\cdot)$, given by $(6.25)$, then it can be established that all the coefficients in $\bar{\xi}_n^*(\cdot)$ are close in probability to that of $\tilde{\xi}_n(\cdot)$, unlike the case of naive studentized bootstrapped estimator. One point we want to make here that the term $\bar{Z}_n^*$ which is present in the expression of $\mathbf{\tilde{H}_n^*}$, unlike the expression of $\mathbf{H_n^*}$, introduces important third order terms which are crucial in getting second order correctness. Therefore, Theorem 4.3 follows.

\subsection{Outline of the proof of Theorem 5.1}
See supplementary material Das and Lahiri (2017).

\section{Conclusion}
Second order results of Perturbation Bootstrap method in regression M-estimation are established. It is shown that the classical way of studentization in perturbation bootstrap setup is not sufficient for correcting the distribution of the regression M-estimator upto second order. This is a general statement corresponding to the fact that the usual studentized perturbation bootstrapped estimator is not capable of correcting the effect of skewness of the error distribution in least square regression. Novel modification is proposed in general setup by properly incorporating the effect of the randomization of the random perturbing quantities in the prevalent studentization factor and is shown as second order correct in both IID and non-IID error setup. Thus, in a way the results in this paper establish perturbation bootstrap method as a refinement of the approximation of the exact distribution of the regression M-estimator over asymptotic normality. The second order result in non-IID case establishes robustness of the perturbation bootstrap towards the presence of heteroscedasticity, similar to the wild bootstrap, but in the more general setup of M-estimation.
This is an important finding from the perspective of S.O.C. inferences regarding the regression parameters.

\section*{Acknowledgement}

The authors would like to thank the two referees, the associate editor and the editor for many constructive comments. They encouraged the authors to add a section on the performance of perturbation bootstrap when the  errors are heteroscedastic (Section 5).

\begin{supplement} 
\stitle{Supplement to ``Second Order Correctness of Perturbation Bootstrap M-Estimator of 
Multiple Linear Regression Parameter''}
\sdatatype{.pdf}
\sdescription{Details of the proofs are provided.}
\end{supplement}

\end{document}